\documentclass[psamsfonts]{amsart}
\usepackage[utf8]{inputenc}
\usepackage{amsfonts,bbm, soul}
\usepackage{hyperref}
\hypersetup{
pdftitle={On The Quasitrace Problem and a Characterization of W*-algebras},
pdfsubject={Mathematics, Operator Algebras, C*-algebras},
pdfauthor={Gow, Alec},
pdfkeywords={Quasitraces, C*-algebras, AW*-algebras, traces, von Neumann algebras, 46L05, 46L35}
}
\usepackage{amsmath}
\usepackage{xcolor}
\usepackage{amsthm}
\usepackage{pdflscape}
\usepackage{mathrsfs}
\usepackage{mathtools}
\usepackage{enumerate}
\usepackage{xcolor}

\usepackage{latexsym,amsfonts}
\usepackage{amsbsy}
\usepackage{amssymb,amsmath}
\usepackage{amscd}
\usepackage{mathrsfs}
\usepackage{physics}
\usepackage{xurl}
\usepackage{lineno}

\hypersetup{
    colorlinks=false,       %
    urlbordercolor={1 1 1}  %
}

\newtheorem{mainthm}{Theorem}

\newtheorem{thm}{Theorem}[section]
\newtheorem{cor}[thm]{Corollary}
\newtheorem{prop}[thm]{Proposition}
\newtheorem{lem}[thm]{Lemma}
\newtheorem{conj}[thm]{Conjecture}
\newtheorem{quest}[thm]{Question}

\theoremstyle{definition}
\newtheorem{defn}[thm]{Definition}

\newtheorem{exmp}[thm]{Example}

\theoremstyle{remark}
\newtheorem{rem}[thm]{Remark}
\newtheorem{rems}[thm]{Remarks}
\newtheorem{warn}[thm]{Warning}

\newcommand{\B}{\mathcal{B}}
\newcommand{\A}{\mathcal{A}}
\newcommand{\Asa}{\mathcal{A}_{sa}}

\DeclareMathOperator{\Proj}{Proj}

\makeatletter
\let\c@equation\c@thm
\makeatother
\numberwithin{equation}{section}

\allowdisplaybreaks

\title[Quasitraces and W*-algebras]{On The Quasitrace Problem and a Characterization of W*-algebras}

\author{Alec Gow}
\address{Department of Pure Mathematics and Institute for Quantum Computing, University of Waterloo, 200 University Ave W, Waterloo, Ontario, Canada, N2L 3G1}
\email{a2gow@uwaterloo.ca}

\begin{document}

\begin{abstract}
We conjecture that a unital C*-algebra is a W*-algebra if and only if each of its maximal abelian self-adjoint subalgebras is a W*-algebra; this is a space-free analogue of a known result due to G.K. Pedersen. Our main result is a proof that this conjecture holds for finite C*-algebras if and only if every $2$-quasitrace on a unital C*-algebra is a trace. We also show that the spatial condition in Pedersen's Theorem can be substantially weakened for AW*-factors. Finally, we give a new characterization of Type II$_1$ W*-factors among Type II$_1$ AW*-factors, which allows us to relate the question of (quasi)linearity of functionals on finite AW*-algebras to the question of monotone completeness of AW*-algebras.  \\

\end{abstract}

\maketitle
\tableofcontents
\section{Introduction}

 In 1951, Kaplansky \cite{Kap51} introduced AW*-algebras as an abstract generalization of W*-algebras (viz. those C*-algebras that can be faithfully represented as von Neumann algebras on some Hilbert space). While Kaplansky's original goal was to give a purely algebraic and intrinsic characterization of W*-algebras, without making reference to Hilbert space representations, it soon became clear that W*-algebras are a proper subclass of all AW*-algebras. In particular, Dixmier \cite{Dixmier} showed that these two classes of C*-algebras do not coincide, even in the commutative case. 
The C*-algebra of continuous functions on a compact Hausdorff space $X$  is an AW*-algebra if and only if $X$ is \emph{Stonean}; it is a W*-algebra if and only if $X$ is \emph{hyperstonean}. %

Despite significant study from the 1950s through the 1980s\footnote{ 
The study of AW*-algebras enjoyed renewed interest in the 2010s (see e.g., \cite{Hamhalter}, \cite{HR1}, \cite{HR2}), thanks to their appearance in the \emph{Bohrification} program of Landsman and others (see \cite{Landsman} and references therein). In this research program, it became apparent that $Comm(\mathcal{A})$, the poset of commutative subalgebras of a C*-algebra $\mathcal{A}$ ordered by set inclusion, was a powerful tool that could be used to study $\A$ itself. It turns out that $Comm(\mathcal{A})$ is a reasonably ``good'' invariant for $\A$ when $\A$ is an AW*-algebra. See \cite{Lindenhovius} for a complete account of the study of $Comm(\mathcal{A})$.}, the theory of AW*-algebras remains less well-understood than that of W*-algebras (and C*-algebras, more generally) which are central objects of study in modern analysis.
One of the major unanswered problems in the theory of AW*-algebras is the question of when an AW*-algebra with W*-centre is itself a W*-algebra. 

A particular case of this question is when exactly an \emph{AW*-factor}, i.e., an AW*-algebra with trivial centre, is a W*-algebra. 
In one of his original papers, Kaplansky  \cite{Kap52} showed that every Type I AW*-factor, indeed every Type I AW*-algebra with W*-centre, is itself a W*-algebra. 
In the 1970s, Dyer \cite{Dyer} and Takenouchi \cite{Takenouchi} independently showed that Type III AW*-factors may not be von Neumann algebras, as they need not admit any normal states.
While it is known that any Type II AW*-factor that admits a faithful state is a W*-algebra (see \cite{Wright75} for the Type II$_1$ case; \cite{CP} for the Type II$_\infty$ case), the question of whether every Type II AW*-factor admits a faithful state is open. 
A longstanding conjecture, originally due to Kaplansky, posits that the answer is yes.
\begin{conj}[Kaplansky] \label{KapConjecture}
Every Type II$_1$ AW*-factor is a W*-factor.
\end{conj}

If Kaplansky's Conjecture has a positive solution, then by standard W*-algebra theory, every Type II$_1$ AW*-factor must admit a unique tracial state, which we will simply call a \emph{trace}: a (normalized) positive linear functional $\tau: \mathcal{A} \to \mathbb{C}$ such that $\tau(x^*x) = \tau (xx^*)$ for all $x \in \mathcal{A}$. 
As it stands, we only know that Type II$_1$ AW*-factors admit more generalized maps called \emph{quasitraces}. 
These are (normalized) maps $\tau: \mathcal{A} \to \mathbb{C}$ that are tracial in the sense that $0 \leq \tau(x^*x) = \tau (xx^*)$ for all $x \in \mathcal{A}$, but are only required to be linear on commutative C*-subalgebras of $\mathcal{A}$ and to satisfy the property that $\tau(a+ib) = \tau(a) + i\tau(b)$ when $a, b$ are self-adjoint elements. If a quasitrace on $\A$ can be extended to a quasitrace on $M_n(\mathcal{A})$ then it is called an \emph{$n$-quasitrace}. 

Through studying which unital C*-algebras admit $2$-quasitraces, we move from the setting of Type II$_1$ AW*-factors to the classification of C*-algebras (see, e.g., \cite{RClassification} for the role of (quasi)traces in classification theory).
By the work of Cuntz \cite{Cuntz} and Handelman \cite{Handelman}, we know that a simple unital C*-algebra admits a $2$-quasitrace if and only if it is \emph{stably finite}. 
But does every stably finite unital C*-algebra admit a tracial state? 
It turns out that a positive answer to this question is equivalent to Kaplansky's Conjecture.
In particular, Blackadar and Handelman \cite{BH} proved that Kaplansky's Conjecture is equivalent to the following.

\begin{conj}[2Q]\label{2Q}
   Every $2$-quasitrace on a unital C*-algebra is a trace. 
\end{conj}

The best partial result towards resolving the 2Q Conjecture was given by Haagerup \cite{Haagerup}. 
By a clever construction, he extended the work of Blackadar and Handelman \cite{BH} to show that every $2$-quasitrace on an \emph{exact} C*-algebra is a trace. 

Recall that Kaplansky's motivation for introducing AW*-algebras was an attempt to give an intrinsic characterization of W*-algebras. 
While he was ultimately unsuccessful in this endeavour, positive results have been established in other work. 
The most famous result along these lines is the elegant theorem of Sakai \cite{Sakai}, which says that a C*-algebra is a W*-algebra if and only if it admits a unique isometric Banach space predual. 
For more recent work in this direction, see \cite{Pham}.

Continuing along this line of inquiry, we consider a natural space-free (i.e., representation-free) analogue of a theorem of G.K. Pedersen \cite{Pedersen} which states that a unital C*-algebra $\mathcal{A} \subseteq \mathcal{B}(H)$ is a concrete von Neumann algebra on $H$ if and only if every maximal abelian self-adjoint subalgebra (MASA) of $\A$ is a concrete von Neumann algebra on $H$.
In particular, we conjecture that a W*-algebra can be completely characterized by its MASAs.  

\begin{conj}[W*-Pedersen]\label{WPC}
   A unital C*-algebra $\A$ is a W*-algebra if and only if every maximal abelian self-adjoint subalgebra of $\A$ is a W*-algebra.  
\end{conj}

In our main result, we establish that Kaplansky's Conjecture and the 2Q Conjecture are equivalent to the W*-Pedersen Conjecture for finite C*-algebras.
In particular, we prove the following result. (The equivalence of $(i)$ and $(ii)$ is known; our contribution is the proof of their equivalence with $(iii)$ and $(iv)$.)

\begin{mainthm}[Theorem \ref{MainThm}]
    The following are equivalent:
    \begin{enumerate}[(i)]
        \item Every Type II$_1$ AW*-factor is a Type II$_1$ W*-factor.
        \item Any 2-quasitrace on a unital C*-algebra is a trace.
        \item A finite C*-algebra is a W*-algebra if and only if all of its MASAs are W*-algebras.
        \item Every Type II$_1$ AW*-factor is a $\ast$-homomorphic retract of its bidual.
    \end{enumerate}
\end{mainthm}

Condition $(iv)$ is established by considering the more general study of quasilinear maps. 
The best-known result in that direction is the generalized Mackey-Gleason Theorem of Bunce and Wright \cite{BunceWright}, which essentially states that any bounded quasilinear map from a W*-algebra (with no Type I$_2$ direct summand) to a Banach space is a bounded linear operator. 
Their work is the key ingredient that allows us to establish that a Type II$_1$ AW*-factor $\A$ is a W*-factor if and only if it is a  \emph{$\ast$-homomorphic retract of its bidual}.

It turns out that studying which unital C*-algebras arise as $\ast$-homomorphic retracts of their biduals is related to the other major open problem in the theory of AW*-algebras. We will say more about this below.

\smallskip
While we were unable to obtain a proof of our conjectured characterization of W*-algebras, we were able to obtain an intermediate result between our abstract conjecture and the concrete Pedersen Theorem for certain simple unital C*-algebras. In particular, we prove the following strengthening of Pedersen's Theorem for simple AW*-factors.

\begin{mainthm}[Theorem \ref{mainTheorem2}]
     Let $\A$ be a simple unital C*-algebra in which every family of mutually orthogonal projections is at most countably infinite. Then $\A$ is a W*-factor if and only if \begin{enumerate}[(i)]
        \item every MASA of $\A$ is an AW*-algebra, and 
        \item there exists a $\ast$-representation $(\pi, H)$ of $\A$ such that for some MASA $M$ of $\mathcal{A}$, $\pi|_M$ is non-zero and completely additive on projections.
    \end{enumerate}
\end{mainthm}

 We obtain a similar result in the case of non-simple unital C*-algebras with trivial centre. 

\begin{mainthm}[Theorem \ref{mainThm3}]
    Let $\A$ be a non-simple, unital  C*-algebra with trivial centre and in which every family of mutually orthogonal projections is at most countably infinite. 
    Then $\A$ is a W*-factor if and only if
    \begin{enumerate}[(i)]
        \item every MASA of $\A$ is an AW*-algebra, and 
        \item there exists a MASA $M$ which contains a non-zero finite projection $e$ and a $\ast$-representation $(\pi, H)$ of $\A$ such that $\pi|_M$ is non-zero, completely additive on projections, and $\pi(e) \neq 0$.
    \end{enumerate}
\end{mainthm}

Our main Pedersen-type results are established using only standard AW*-algebra theory. 
Much of the previous work along these lines used the more general theory of \emph{AW*-modules} -- introduced by Kaplansky in \cite{Kaplansky3} -- to describe when an AW*-algebra is a W*-algebra (e.g. \cite{Berberian83-1}, \cite{Berberian83-2}, \cite{ESW}).    

\smallskip 
Finally, as mentioned above, we briefly consider the other major open problem in the theory of AW*-algebras. 
In the case of commutative C*-algebras, it is well-known that being an AW*-algebra is the same as being \emph{monotone complete}. 
Moreover, every monotone complete C*-algebra is an AW*-algebra. 
The converse is an open problem.

\begin{conj}[Monotone Completeness]
    Every AW*-algebra is monotone complete.
\end{conj}

In the final section of the paper, motivated by a MathOverflow question raised by Thiel and a subsequent answer by Gao \cite{MO}, we establish a bridge between the two major open problems in the theory of AW*-algebras.
In particular, we show that if every monotone complete C*-algebra is a retract of its bidual, then via condition (iv) of Theorem A, the problem of resolving Kaplansky's Conjecture (and the equivalent conjectures) reduces to a subcase of the problem of resolving the Monotone Completeness Conjecture.  
\smallskip

In order to make this paper accessible to operator algebraists without special knowledge of AW*-algebras or quasitraces, we have included significant background material. Experts in the former may wish to skip Section \ref{subsection:AW}; experts in the latter may wish to skip Section \ref{subsection:dimension} and Section \ref{HC}.
We have also added details to proofs of some previously known results in order to improve clarity, harmonize terminology, and keep this paper as self-contained as possible. 

 The paper is organized as follows. We begin in Section \ref{Background} with a review of the necessary mathematical background on AW*-algebras, quasitraces, and Haagerup's partial solution to Kaplansky's Conjecture. 
 In Section \ref{PedersenSection}, we consider the general problem of determining the structure of a C*-algebra by its MASAs, give a detailed proof of Pedersen's Theorem, and provide some evidence that the space-free W*-Pedersen Conjecture may have a positive solution.
 In Section \ref{MainSection}, we prove our first main theorem. We prove the first new equivalence by proving that the Finite W*-Pedersen Conjecture implies the Kaplansky Conjecture, sketching the known result that the latter implies the 2Q Conjecture, and proving that the 2Q Conjecture implies the Finite W*-Pedersen Conjecture. We then directly prove that the bidual-retract condition is equivalent to Kaplansky's Conjecture.
 In Section \ref{ProperlyInfiniteSec}, we show that a Mackey-Gleason result holds for properly infinite AW*-algebras and we obtain the strengthened Pedersen-type theorems for countably decomposable AW*-factors. 
 In Section \ref{MCQuestion}, we establish a path by which we could view Kaplansky's Conjecture as a specific case of the Monotone Completeness Conjecture. 

 \smallskip 
\noindent \textbf{Acknowledgements.} The author thanks Ken Davidson for his comments on an earlier version of this paper; Hannes Thiel for suggesting that the gap identified in the author's previous work could be equivalent to the central problem studied here, and for his helpful comments on this revision; and Michael Brannan for many helpful conversations, edits, and his supervision. 

\section{Background} \label{Background}

We begin with a review of well-known facts about our key ingredients -- AW*-algebras, quasitraces, and Haagerup's construction -- providing references to the relevant literature throughout.

In the sequel, $\A$ will always denote a unital C*-algebra. We write $\Asa$ for the partially ordered, real Banach space of self-adjoint elements of $\mathcal{A}$, $\Proj(\mathcal{A})$ for the set of projections (i.e., self-adjoint idempotents) of $\mathcal{A}$, $\mathcal{A}_1$ for the unit ball of $\mathcal{A}$, and $\mathcal{A}_+$ for the positive elements of $\mathcal{A}$. If $\mathcal{A} = \mathcal{B}(H)$, the bounded operators on some Hilbert space $H$, then we denote the set of projections by $\mathcal{P}(H)$. Finally, whenever we speak of a representation of a C*-algebra, we shall always mean a non-degenerate $\ast$-representation. We will simply call a representation of a C*-algebra on a separable Hilbert space a separable representation; if a C*-algebra admits a faithful separable representation, we will say that it is \emph{separably representable}.

\subsection{AW*-algebras}\label{subsection:AW}

In addition to the original papers of Kaplansky (\cite{Kap51}, \cite{Kap52}), the canonical reference for AW*-algebras is the textbook of Berberian \cite{Berberian}. We begin with Kaplansky's original definition.

\begin{defn}
   An \emph{AW*-algebra} $\A$ is a unital C*-algebra such that
   \begin{enumerate}[(i)]
       \item any set of pairwise orthogonal projections has a least upper bound (LUB), and
       \item any maximal abelian self-adjoint subalgebra (MASA) is generated as a subalgebra by its projections.
   \end{enumerate}
\end{defn}

\begin{exmp}\cite[Theorem 2.4]{Kap51}
    If $\A$ is an AW*-algebra, then so is its centre ${\mathcal{Z(A)}:= \{z \in \mathcal{A} : za = az, \hspace{0.01cm} \forall  a \in \mathcal{A}\}}$, any MASA of $\mathcal{A}$, and any corner subalgebra $p\mathcal{A}p$ where $p \in \Proj({\mathcal{A}})$.
\end{exmp}

There is another characterization of AW*-algebras which is of particular interest to us, as it has the same ``flavour'' as Pedersen's Theorem \ref{PTheorem}. 
This characterization was a longstanding folklore result in the AW*-algebra literature; a proof was first published by Sait{\^o} and Wright in 2015 \cite{SW}.

\begin{prop}\label{AWcharacterization}\cite[Proposition 1.4]{SW}
    A unital C*-algebra $\A$ is an AW*-algebra if and only if every MASA of $\A$ is an AW*-algebra.
\end{prop}

\begin{rem}
    The original statement of Sait{\^o} and Wright's characterization says that $\A$ is an AW*-algebra if and only if every MASA is \emph{monotone complete}, i.e., if every norm-bounded monotone increasing net in $\mathcal{A}_{sa}$ has a supremum in $\mathcal{A}_{sa}$.
    Every monotone complete C*-algebra is an AW*-algebra but the converse is an open question. 
    (Indeed, it is \emph{the} major open problem in the theory of AW*-algebras, apart from the resolution of Kaplansky's Conjecture. We consider the relationship between the two questions in Section \ref{MCQuestion}.) 
    It is known, however, that the notions of monotone completeness and AW*-ness coincide for commutative (unital) C*-algebras. 

    In particular, a commutative C*-algebra $\mathcal{C}(X)$ is monotone complete (equivalently, an AW*-algebra) if and only if the compact Hausdorff space is \emph{extremally disconnected}, i.e., the closure of any open subset of $X$ is again open. 
    For compact Hausdorff spaces, extremally disconnected is also called \emph{Stonean} in the literature. 
    We will use the latter, as it contrasts nicely with the case for W*-algebras; recall that $\mathcal{C}(X)$ is a W*-algebra if and only if $X$ is \emph{hyperstonean} (i.e., $X$ is Stonean and admits sufficiently many normal Radon measures).
\end{rem}

We will have occasion to use the fact that when $\A$ is an AW*-algebra, $\Proj(\mathcal{A})$ is a complete lattice (see, e.g., Corollary 3, \cite{Kap51}), meaning that we can take arbitrary suprema and infima of families of (not necessarily mutually orthogonal) projections. Given a family of projections $(p_i)_{i \in I}$ in $\mathcal{A}$ with LUB $p$ and greatest lower bound (GLB) $q$, we write 
\[ p = \bigvee p_i = \sup_{\Proj({\mathcal{A}})} \{p_i\}, \hspace{0.25cm} q = \bigwedge p_i = \inf_{\Proj ({\mathcal{A}})} \{p_i\}. \]
We will often drop the $\Proj({\mathcal{A}})$ subscript notation if it is clear from context that we are taking suprema/infima of projections in $\mathcal{A}$.

If $\A$ is a von Neumann algebra and $(p_i)_{i \in I}$ a family of pairwise orthogonal projections in $\mathcal{A}$, then $\sum_{i} p_i$ is well-defined as the SOT-limit of finite sums (see e.g., \cite[Corollary 2.7.2]{Peterson}). Moreover, the SOT-limit coincides with the LUB of $\{ p_i\}$ in $\Proj(\A)$, \[p = \sum_{i} p_i =  \bigvee p_i = \sup_{\Proj({\mathcal{A}})} \{p_i\}.\]

Much of the usual terminology from von Neumann algebras carries over to AW*-algebras. 
For example, we say that an AW*-algebra is a \emph{factor} if it has trivial centre, i.e.,  $\mathcal{Z}(\mathcal{A}) = \mathbb{C}1$.
In order to make things clear, if $\A$ is an AW*-algebra which is also a factor, we will call it an AW*-factor. If $\A$ is a W*-algebra which is also a factor, we will call it a W*-factor.

When $\A$ is an AW*-algebra, the fact that $\Proj(\mathcal{A})$ is a complete lattice means that for any projection $p \in \Proj(\mathcal{A})$, the set $$\{z \in \Proj(\mathcal{Z}(\mathcal{A})): z \geq p\}$$ has an infimum in $\mathcal{A}$, which we call the \emph{central cover} of $p$, denoted $z(p)$ or $z_{\mathcal{A}}(p)$, where the subscript emphasizes the particular AW*-algebra in which the central cover is computed. (For a direct proof that $z(p)$ is indeed a central projection, see \cite[Lemma 8.6.7]{Lindenhovius}. Alternatively, Example \ref{cover} below sheds light on why this is the case.)

Like their von Neumann algebra counterparts, AW*-algebras have a type decomposition, which is largely based on properties of projections. 
The same notion of equivalence of projections carries over to this setting.
Recall that two projections $p,q$ are said to be \emph{(Murray-von Neumann) equivalent}, written $p\sim q$, if there is a partial isometry $u$ such that $p = uu^*$ and $q = u^*u$. 
We say that a projection $p$ is \emph{subequivalent} to another projection $q$ if there exists another projection $q_0$ such that $p \sim q_0 \leq q$. In this case, we write $p \lesssim q.$
We say that a projection $p$ is \emph{finite} if $p \sim q \leq p$ implies $p = q$; otherwise, we say that $p$ is \emph{infinite}. 

\begin{defn}
    An AW*-algebra $\A$ is said to be
    \begin{enumerate}[(i)]
        \item \emph{discrete} if there exists $p \in \Proj(\mathcal{A})$ such that $z(p) = 1$ and $p \mathcal{A} p$ is abelian. 
        \item \emph{continuous} if $p = 0$ is the only element of $\Proj(\mathcal{A})$ such that $p\mathcal{A}p$ is abelian. 
        \item \emph{finite} if $1 \in \Proj(\mathcal{A})$ is finite. That is, whenever $x \in \mathcal{A}$ satisfies  $x^*x = 1$ then $xx^* = 1$. 
        \item \emph{semifinite} if there exists $p \in \Proj(\mathcal{A})$ such that $p\mathcal{A}p$ is finite and $z(p) = 1$. 
        \item \emph{properly infinite} if every $0 \neq p \in \Proj(\mathcal{Z(A)})$ is infinite.
        \item \emph{purely infinite} if $p = 0$ is the only element of $\Proj(\mathcal{A})$ such that $p\mathcal{A}p$ is finite. 
    \end{enumerate}
\end{defn}

\noindent Notice that the notion of finiteness (in the sense that $x^*x = 1 \iff xx^*=1$) is well-defined for any unital C*-algebra. 

\begin{defn}
    An AW*-algebra $\A$ is said to be
    \begin{enumerate}[(i)]
    \item \emph{Type I$_{\text{fin}}$} if it is discrete and finite.
    \item \emph{Type I$_{\infty}$} if it is discrete and properly infinite.
    \item \emph{Type II$_1$} if it is continuous and finite. 
    \item \emph{Type II$_\infty$} if it is continuous, semifinite, and properly infinite. 
    \item \emph{Type III} if it is purely infinite.
    \end{enumerate}   
\end{defn}

\noindent Any AW*-algebra $\mathcal{A}$ can be written as the direct sum of parts of each type. That is, $\mathcal{A} = \mathcal{A}_{I_{\text{fin}}} \oplus \mathcal{A}_{I_{\infty}} \oplus \mathcal{A}_{II_1} \oplus \mathcal{A}_{II_\infty} \oplus \mathcal{A}_{III}$, where $\mathcal{A}_{\nu}$ is an AW*-algebra of Type $\nu$ and $\nu \in \{\, \text{I$_{\text{fin}}$, I$_{\infty}$, II$_{1}$, II$_{\infty}$, III}\,\}$. 

One can decompose the Type I$_\text{fin}$ AW*-algebras further, by showing that they are direct sums of Type I$_n$ AW*-algebras, where $n \in \mathbb{N}.$ 
Kaplansky showed \cite{Kap52} that if $\A$ is an AW*-algebra of Type I$_n$ then it is isomorphic to $M_n(\mathcal{Z}(\A))$. 
In particular, a Type I$_n$ AW*-factor is of the form $M_n(\mathbb{C})$.

For AW*-algebras, being purely infinite is a strictly stronger property than being properly infinite. 
When we refer to a \emph{properly infinite AW*-algebra}, we mean an AW*-algebra that can only have direct summands of Type I$_\infty$, Type II$_\infty$, and/or Type III.  
It is useful to know that, in any properly infinite AW*-algebra, there exists $p \in \Proj(\mathcal{A})$ such that $1 \sim p \sim 1-p$. 

We refer to any AW*-algebra without a Type III direct summand as a \emph{semifinite AW*-algebra.}

\begin{rem}\label{rem:infiniteTerminology}
    The terminology used to describe infinite AW*-algebras varies in the older literature. 
    For example, what Kaplansky defines as ``purely infinite'' in his original paper \cite{Kap51} refers to the condition that we, following modern conventions, call ``properly infinite''. 
\end{rem}

As in the von Neumann algebra case, all finite factors are simple C*-algebras. 
The Type III factors that are \emph{countably decomposable}, meaning those in which any family of pairwise orthogonal projections is at most countably infinite, are also simple. 
On the other hand, AW*-factors of Type I$_\infty$ or Type II$_\infty$ are not simple. 
We note that for von Neumann algebras, simplicity for Type III factors is equivalent to countable decomposability, but there exist (wild) simple Type III AW*-factors which are not countably decomposable (see e.g., \cite[Lemma 5.2]{Hamana}).

\smallskip
An important concept throughout this paper will be the idea of compatibility of the AW*- (and/or W*-) structures of C*-algebras. 
In general, if $\mathcal{A}, \mathcal{B}$ are both AW*-algebras and $1 \in \mathcal{B} \subseteq \mathcal{A}$ is a unital inclusion of C*-algebras, the AW*-structure of $\B$ need not be the same as the AW*-structure of $\mathcal{A}$. 

\begin{exmp}\label{exmp:universalW}
    Let $\B$ be a W*-algebra that is concretely represented on its universal Hilbert space $H_u$. 
    It is well-known that if $\B$ is not finite-dimensional, then $\B$ is not a von Neumann algebra on $H_u$.
    In particular, as $\B$ need not be SOT-closed in $\mathcal{B}(H_u)$, there exists a family of orthogonal projections $(p_i)_{i \in I}$ in $\B$ such that $\sup_{\Proj({\mathcal{B}})} \{p_i\} \neq \sup_{\mathcal{P}(H_u)}\{p_i\}$. 
\end{exmp}

\noindent In order for the AW*-structure(s) of $\A$ and $\B$ to be fully compatible, we need one to be an \emph{AW*-subalgebra} of the other.

\begin{defn}
    Let $\A$ be an AW*-algebra. We say that a C*-subalgebra $\mathcal{B} \subseteq \mathcal{A}$ is an \emph{AW*-subalgebra} of $\A$ if 
    \begin{enumerate}[(i)]
        \item $\mathcal{B}$ is an AW*-algebra, and 
        \item given any family of mutually orthogonal projections $(p_i)_{i \in I} \subseteq \mathcal{B},$ the LUB of $(p_i)$ in $\Proj({\mathcal{B}})$ is equal to the LUB of $(p_i)$ calculated in $\Proj ({\mathcal{A}})$, i.e., \[ \sup_{\Proj ({\mathcal{B}})} \{p_i\} = \sup_{\Proj ({\mathcal{A}})} \{p_i\}.\] 
    \end{enumerate}
\end{defn}

In fact, compatibility of suprema described in (ii) holds for arbitrary families of projections.
This is a special case of the more general fact that $\ast$-homomorphisms between AW*-algebras which preserve suprema of families of orthogonal projections preserve suprema of arbitrary families of projections (see e.g., \cite[Lemma 8.1]{HR1}).

\begin{exmp}\label{cover}
    For any AW*-algebra $\mathcal{A}$, $\mathcal{Z}(\mathcal{A})$ is an AW*-subalgebra of $\mathcal{A}$. This means that for any projection $p \in \Proj(\mathcal{A})$, the central cover of $p$ can be computed as the infimum of the set $\{z \in \Proj(\mathcal{Z}(\mathcal{A})): z \geq p\}$ in either AW*-algebra.
\end{exmp}

\begin{exmp}\cite[Lemma 1.6]{SW}\label{exmp:MASA}
    Let $\A$ be an AW*-algebra and let $M$ be a MASA of $\A$. Then $M$ is an AW*-subalgebra of $\A.$
\end{exmp}

\begin{rem} 
The property of ``being an AW*-subalgebra'' is transitive: if $\mathcal{C}$ is an AW*-subalgebra of $\B$ and $\B$ is an AW*-subalgebra of $\mathcal{A}$, then $\mathcal{C}$ is an AW*-subalgebra of $\mathcal{A}$.
\end{rem}

Finally, we adopt the terminology of \cite{PMonograph} in order to describe an AW*-algebra satisfying some nice properties. 

\begin{defn}
     Let $H$ be a Hilbert space and let $\mathcal{A} \subseteq \mathcal{B}(H)$ be an AW*-algebra. We say that $\A$ is a \emph{concrete AW*-algebra} on $H$ if 
     \begin{enumerate}[(i)]
         \item $\A$ contains the spectral projections of each element $a \in \mathcal{A}_{sa}$, and 
        \item for any pairwise orthogonal family of projections $(p_i)_{i \in I}$ in $\mathcal{A}$, $\sum_i p_i \in \mathcal{A}$, where $\sum_i p_i$ is calculated as the SOT-limit of finite sums in $\mathcal{B}(H)$.
     \end{enumerate}
\end{defn}

It is immediate that a concrete AW*-algebra $\mathcal{A} \subseteq \mathcal{B}(H)$ is an AW*-subalgebra of $\mathcal{B}(H)$. 
In fact, these notions coincide. To show this, we use a result due to Feldman \cite{Feldman}. Where Feldman used the (now outdated) notion of \emph{AW*-embedding}, we have updated the terminology to use the equivalent notion of AW*-subalgebra and added more detail to the proof, for clarity.

\begin{lem}\cite[Lemma 1]{Feldman}\label{lem:spectralProj}
    If $\A$ is an AW*-subalgebra of $\mathcal{B}(H)$, then $\A$ contains all the spectral projections of its self-adjoint elements.
\end{lem}

\begin{proof}
    First, suppose that $\A$ is a commutative AW*-algebra, so that $\iota(\mathcal{A}) = \mathcal{C}(X)$ where $X$ is Stonean (that is, compact, Hausdorff, and extremally disconnected) and $\iota: \mathcal{A} \to \mathcal{C}(X)$ is the canonical isomorphism.
    For $a \in \mathcal{A}_{sa}$, let $X(a,r)$ denote the interior of the set $\{x \in X: \iota({a})(x) \leq r\}$.
    Since $X$ is extremally disconnected, $X(a,r)$ is clopen and so its characteristic function $\mathbf{1}_{X(a,r)}$ is continuous.
    Now, let $E_a(r) := \iota^{-1}(\mathbf{1}_{X(a,r)})$.
    Then $E_a(r)$ is a monotone increasing, projection-valued function of $r$ with $\sup E_a(r) = 1$ and $\inf E_a(r) = 0$.
   
    Let $p \in \Proj(\mathcal{A})$ be a projection such that $p \leq E_a(r)$ for all $r > r_0$. 
    Then $\iota(p)$ is the characteristic function of some clopen set $X_0$, where $X_0 \subseteq X(a,r)$ for every $r > r_0$.
    This means that for each $x \in X_0$, we have $\iota(a)(x) \leq r_0$ and so by definition $x \in X(a, r_0)$ and therefore $p \leq E_a(r_0).$
    So, $E_a(r_0)$ is the infimum of $E_a(r)$ when $r > r_0$, meaning that $E_a(r)$ (as a function of $r$) is order-continuous from above.

    Now, let $m := \inf_X \iota(a)$ and let $M := \sup_X \iota(a)$.
    If $m = M$, then $\iota(a) = m \cdot 1_{\mathcal{C}(X)}$, hence $a = m \cdot 1_\A$ and there is nothing to prove. 
    Otherwise, fix some $r_0 < m$ and let us also fix a partition $m = r_1 < r_2 < \cdots < r_{n+1} = M$ of the range of $\iota(a)$. 
    Define finite sums \[S_n :=  \sum_{k=0}^{n} r_k \big(E_a(r_{k+1}) - E_a(r_k)\big).\]
    We define pairwise disjoint, clopen sets  
    \[ A_k := X(a,r_{k+1}) \setminus X(a,r_k), \quad 0 \le k \le n.\] 
    Since $r_0 < m$, we have $X(a, r_0) = \emptyset$ and so $A_0 = X(a, m) = X(a, r_1)$.
    Moreover, $X(a, r_{n+1}) = X(a, M) = X$ so we have
    \[\bigcup_{k=0}^n A_k = X. \]

    Let $s_n(x) := \iota(S_n)(x).$
    For $x \in A_k$ we see that 
    \begin{align*}
        s_n(x) = \iota(S_n)(x) &= \sum_{\ell=0}^{n} r_\ell \big(\mathbf{1}_{X(a,r_{\ell+1})} - \mathbf{1}_{X(a,r_\ell)}\big)(x) \\
        &= r_k.
    \end{align*}

    Whenever $x \in X(a,r_{k+1})$, there is a neighbourhood of $x$ on which $\iota(a) \le r_{k+1}$, hence $\iota(a)(x) \le r_{k+1}$.  
    If $x \notin X(a,r_k)$, then $x$ lies in the closure of the open set $\{y \in X: \iota(a)(y) > r_k\}$, so there must be some net $(x_i)_{i \in I} \subseteq X$ such that $x_i \to x$ and $\iota(a)(x_i) > r_k$. 
    By continuity, $\iota(a)(x) \geq r_k$ and so for all $x \in A_k$, we have
    \begin{align*}
        0 \leq \iota(a)(x) - r_k &\leq r_{k+1} - r_k \\
        \iff 0 \leq (\iota(a) - r_k)\mathbf{1}_{A_k} &\leq (r_{k+1}-r_k)\mathbf{1}_{A_k}.
    \end{align*}
    Letting $\delta = \max\limits_{0 \le k \le n}(r_{k+1}-r_k)$ and summing over $k$, we get
    \begin{align*}
        0 \leq \iota(a) - s_n &\leq \delta\cdot \mathbf{1}_X \\
        \implies \norm{a-S_n} = \norm{\iota(a) - s_n}_\infty &\leq \delta.
    \end{align*}
    Taking $r_0$ arbitrarily close to $m$ and $\delta \to 0$, we get that $S_n$ approximates $a$ uniformly.
    Therefore, 
    \[ a = \int_{m}^M r \,dE_a(r) = \int_\mathbb{R} r \,dE_a(r),\] as a uniform Riemann-Stieltjes integral.

    Now, suppose that $\A$ is any AW*-algebra which is an AW*-subalgebra of $\mathcal{B}(H)$.
    For any $a \in \mathcal{A}_{sa},$ there is some MASA $M$ of $\mathcal{A}$ that contains $a$.
    So, we can find a family of projections $E_a( \cdot) \subseteq M$ as above.
    Since $M \subseteq \mathcal{A} \subseteq \mathcal{B}(H)$, where each inclusion respects the AW*-structure, $M$ is also an AW*-subalgebra of $\mathcal{B}(H)$.
    This means that taking suprema (and infima) of projections in $M$ coincides with taking SOT-limits in $\mathcal{B}(H)$. 
    Moreover, since the partial sums $S_n$ uniformly approximate $a$, we get that $E_a(\cdot)$ is the usual spectral resolution of $a$ in $\mathcal{B}(H)$.
    Therefore, $\mathcal{A}$ contains all the spectral projections $E_a(r)$ of any self-adjoint $a \in \mathcal{A}_{sa}$. 
\end{proof}

\begin{rem}
    There is always a well-defined notion of spectral projections for self-adjoint elements in an abstract AW*-algebra $\A$, essentially by the construction outlined in the first part of the proof of Lemma \ref{lem:spectralProj}.
    However, even if $\A$ is faithfully represented as a C*-algebra on $H$, the spectral projections in $\A$ need not coincide with the usual spectral projections in $\mathcal{B}(H)$. 
   The lemma only says that they will coincide if $\A$ is an AW*-subalgebra of $\B(H).$   
\end{rem}

\begin{cor}\label{cor:dominateProjection}
      Let $\A$ be an AW*-algebra. 
      For any $0 \neq a \in \A_+$, there exists a non-zero projection $p \in \Proj(\A)$ and a scalar $\varepsilon > 0$ such that $a \geq \varepsilon p.$  
\end{cor}

\begin{proof}
    Given $0 \neq a \in A_+$, let us denote $\mathcal{C}(X) \cong \iota(M),$ where $M$ is the MASA of $A$ containing $a$ and $\iota: \mathcal{M} \to \mathcal{C}(X)$ is the canonical isomorphism.
    Pick some $\varepsilon \in (0, \norm{\iota(a)})$.
    Let \[ U := \{x \in X: \iota(a)(x) > \varepsilon\}.\]
    Since $X$ is Stonean, $\overline{U}$ is clopen, hence  $p := \iota^{-1}(\mathbf{1}_{\overline{U}})$ is a projection.
    By continuity, we see that $a \geq \varepsilon p.$
\end{proof}

We obtain the following result from the relevant definitions and Lemma \ref{lem:spectralProj}.

\begin{cor}\label{cor:concreteAW}
    Let $\mathcal{A} \subseteq \mathcal{B}(H)$ be an AW*-algebra. Then $\A$ is a concrete AW*-algebra on $H$ if and only if it is an AW*-subalgebra of $\mathcal{B}(H)$.
\end{cor}

Before moving on, we shall describe an example which highlights the delicate interplay between the structure of an abstract W*-algebra $\mathcal{M}$ and the structure of $\pi(\mathcal{M})$, where $(\pi, H)$ is an arbitrary representation of $\mathcal{M}$ as a C*-algebra. 

It is well-known that if $\mathcal{M}$ is a W*-algebra, then it can admit (many) representations $(\pi, H)$ such that $\pi(\mathcal{M})$ is not a concrete von Neumann algebra on $\mathcal{B}(H)$ (as in e.g., Example \ref{exmp:universalW}).
Moreover, it is known that if $(\pi, H)$ is a separable representation of $\mathcal{M}$, then $\pi(\mathcal{M})$ is always an AW*-algebra (by \cite[Proposition 2.3]{Handelman}). 
However, we will show that $\pi(\mathcal{M})$ can still lose its abstract W*-structure. 

\begin{exmp}[A non-W*-representation of a W*-algebra on a separable Hilbert space]\label{exmp:nonWseparable}
    Let $\mathcal{C}(X)$ be the so-called \emph{Dixmier algebra} \cite{Dixmier}.
    That is, take \[\mathcal{C}(X) \cong B([0,1])/M([0,1]),\] where $B([0,1])$ is the set of bounded (complex-valued) Borel functions on $[0,1]$ and $M([0,1])$ is the ideal of functions in $B([0,1])$ with meagre support. 
    This is the canonical example of a commutative AW*-algebra which is not a W*-algebra; in topological terms, recall that this means $X$ is Stonean but not hyperstonean.
    
    It is well-known (see e.g., \cite[III.1.8.5]{Blackadar}) that the Dixmier algebra is separably representable; let $(\sigma, H)$ denote a faithful separable representation of $\mathcal{C}(X)$.
    
    Let $S$ be the discretization of $X$. That is, let $S$ be the space obtained by equipping the underlying set of $X$ with the discrete topology.
    Then by basic general topology, the identity map $S \to X$ extends to a continuous surjection $p: \beta S \to X$, where $\beta S$ is the Stone-\v{C}ech compactification.

    Since $X$ is extremally disconnected, Gleason projectivity (\cite[Theorem 2.5]{Gleason}) gives a \emph{continuous section} \[s: X \to \beta S\] with $p \circ s = \operatorname{id}_X.$
    By duality, this gives rise to a $\ast$-homomorphic \emph{retraction} \[r: \mathcal{C}(\beta S) \to \mathcal{C}(X)\] where $r(f) = f \circ s$.

    Set $M:= \mathcal{C}(\beta S) = \ell^{\infty}(S),$ which is a commutative W*-algebra. Then clearly $(\pi = \sigma \circ r, H)$ is a representation of $M$ on a separable Hilbert space, but \[\pi(M) := \sigma(\mathcal{C}(X)) \cong \mathcal{C}(X),\] which is an AW*-algebra but not a W*-algebra.
\end{exmp}

\begin{rem}
    The preceding example may be known to experts familiar with the notions of projectivity and injectivity in commutative C*-algebras and operator systems.
    We have included it here because it does not seem to have previously appeared in this context. 
    Indeed, in the errata to \cite{Blackadar}, Blackadar writes that ``it appears to be unknown'' whether a representation of a W*-algebra on a separable Hilbert space must remain a W*-algebra; our example confirms that this question has a negative solution.
\end{rem}

\subsection{Dimension Functions and Quasitraces}\label{subsection:dimension}

The central sources on traces and quasitraces in this context are Blackadar and Handelman \cite{BH} and Haagerup \cite{Haagerup}. See also \cite{MilhojRordam} for a recent state-of-the-art survey. 

\begin{defn}
    A \emph{1-quasitrace} $\tau$ on $\A$ is a function $\tau: \mathcal{A}\to\mathbb{C}$ such that
    \begin{enumerate}[(i)]
        \item $0 \leq \tau(x^*x) = \tau(xx^*)$ for all $x \in$ $\mathcal{A}$. 
        \item $\tau(a+ib) = \tau(a) + i\tau(b)$ for all $a, b \in$ $\mathcal{A}_{sa}$. 
        \item $\tau$ is linear on every abelian $C^*$-subalgebra of $\mathcal{A}$. 
    \end{enumerate}
    If there exists a 1-quasitrace $\tau_n: M_n(\mathcal{A}) \to \mathbb{C}$ such that $\tau_n(a \otimes e_{11}) = \tau(a)$ for all $a \in \mathcal{A}$, then $\tau$ is said to be an \emph{n-quasitrace}.
\end{defn}

Blackadar and Handelman showed \cite[Proposition II.4.1]{BH} that any $2$-quasitrace extends to an $n$-quasitrace. In light of this result, we use the term \emph{quasitrace} to refer to any normalized $2$-quasitrace, that is, a $2$-quasitrace with $\tau(1) = 1$.  

Any bounded linear functional $\tau: \mathcal{A} \to \mathbb{C}$ is said to be \emph{tracial} if $\tau(x^*x) = \tau(x x^*)$ for all $x \in \mathcal{A}$. We will call a tracial state a \emph{trace}; in other words, a trace is a quasitrace that is linear on the whole C*-algebra. 

When $\A$ is a unital C*-algebra, we write $QT(\mathcal{A})$ to denote the set of all normalized $2$-quasitraces on $\mathcal{A}$. It is well-known that, provided $QT(\mathcal{A})$ is non-empty, it is a Choquet simplex (see e.g., \cite[Theorem II.4.4]{BH}).

The canonical quasitrace of a finite AW*-algebra is built from its unique \emph{centre-valued dimension function}. We summarize the important properties of this function below; proofs can be found in \cite[Chapter 6]{Berberian}.

\begin{thm}
    Let $\mathcal{A}$ be a finite AW*-algebra. There exists a unique centre-valued dimension function ${D: \Proj(\mathcal{A}) \to \mathcal{Z(A)}}$ satisfying the following properties:

\begin{enumerate}[(i)]
    \item $p \sim q$ if and only if $D(p) = D(q)$.
    \item If $p \perp q$ then $D(p+q) = D(p) + D(q)$.
    \item If $z$ is a central projection, then $D(zq) = zD(q)$ for all $q \in \Proj(\mathcal{A})$.
    \item $p \lesssim q$ if and only if $D(p) \leq D(q)$.
    \item $0 \leq D(p) \leq 1$, $D(1) = 1$, and $D(p) = 0$ if and only if $p = 0$.
    \item If $(p_i)_{i \in I}$ is a family of mutually orthogonal projections, then \begin{align*} 
    D\Bigg(\sup_{\Proj (\mathcal{A})}\{p_i\}\Bigg) &= \sum_{i \in I} D(p_i) = \sup \Bigg\{ \sum_{i \in J} D(p_i): J \subseteq I, J \text{ finite }\Bigg\}.
    \end{align*}
\end{enumerate}
\end{thm}

Using the theory of dimension functions, Blackadar and Handelman \cite[Section II]{BH} introduced by name\footnote{The \emph{notion} of such maps goes back to Murray and von Neumann \cite{MvN}.}  the unique \emph{centre-valued quasitrace} $\tau$ on a finite AW*-algebra, $\tau: \mathcal{A} \to \mathcal{Z(A)}$. Again, we summarize the important properties; proofs are scattered throughout the literature but were collected in \cite{Fehlker}.

\begin{thm} \cite[Theorem 1.27]{Fehlker}
    Let $\A$ be a finite AW*-algebra. There exists a unique function $\tau: \mathcal{A} \to \mathcal{Z(A)}$ satisfying the following properties:
    
    \begin{enumerate}[(i)]
    \item $\tau(x+y) = \tau(x) + \tau(y)$ for any commuting, normal elements $x, y \in \mathcal{A}$.
    \item $\tau(\lambda x) = \lambda \tau(x)$ for all $x \in \mathcal{A}$, $\lambda \in \mathbb{C}$.
    \item $0 \leq \tau(x^*x) = \tau(xx^*)$ for all $x \in \mathcal{A}$.
    \item $\tau(x^*x) = 0$ if and only if $x = 0$. 
    \item $\tau(x+iy) = \tau(x) + i\tau(y)$ for all $x,y \in \mathcal{A}_{sa}$. 
    \item $\tau(z) = z$ for all $z \in \mathcal{Z(A)}$.
    \item $x \leq y $ implies $\tau(x) \leq \tau(y)$ for all $x, y \in \mathcal{A}_{sa}$. 
    \item $\tau$ is norm-continuous.
    \item $\tau |_{\Proj (\mathcal{A})} = D$, the centre-valued dimension function.
\end{enumerate}
\end{thm}

It is clear from the properties above that if $\A$ is a finite AW*-factor, then the (unique) ``centre-valued quasitrace'' on $\A$ is a $1$-quasitrace. By a result of Blackadar and Handelman  \cite[Corollary II.1.10]{BH}, every $1$-quasitrace on an AW*-algebra is a $2$-quasitrace (hence an $n$-quasitrace), so the centre-valued quasitrace is indeed a quasitrace in our sense, as well. Moreover, property (iv) tells us that $\tau$ is \emph{faithful}. 

\begin{defn} \label{normal}
    A quasitrace on an AW*-algebra $\A$ is said to be \emph{normal} if for any pairwise orthogonal set of projections $(p_i)_{i \in I} \subseteq \Proj(\mathcal{A})$, \begin{equation*}
        \tau\Bigg(\sup_{\Proj (\mathcal{A})} \{p_i\}\Bigg) = \sum_{i \in I} \tau(p_i) = \sup \Bigg\{ \sum_{i \in J}\tau(p_i): J \subseteq I, J \text{ finite }\Bigg\}.
    \end{equation*}
\end{defn}

\begin{rem}
    Here we have adopted the term ``normal'' quasitrace, as is standard in the quasitrace literature. 
However, we note that if $\varphi: \mathcal{A} \to \mathbb{C}$ is a positive, bounded linear functional satisfying the condition given in Definition \ref{normal}, it is more common in the AW*-algebraic literature to call this functional \emph{completely additive on projections} (CAP). 

The usual definition of a \emph{normal linear functional} is a positive, bounded linear functional preserving directed limits of bounded nets of self-adjoint elements. 
That is, if $(x_i)_{i \in I}$ is a norm-bounded, monotone increasing net in $\mathcal{A}_{sa}$, then $\varphi$ is normal whenever $\sup \varphi(x_i) = \varphi (\sup \{x_i\})$. 
However, these two notions are closely related.
\end{rem}

It is well-known that if $\A$ is a W*-algebra, then a functional $\varphi: \mathcal{A} \to \mathbb{C}$ is normal (in the usual sense) exactly when it is CAP. In fact, this holds for any monotone complete C*-algebra (see e.g., \cite[Proposition 1, $\S$9.34]{Zsido}). In particular, we have the following result.

\begin{lem}\label{normalCAP}
    Let $\A$ be an AW*-algebra and $\varphi: \mathcal{A} \to \mathbb{C}$ a CAP functional. Then $\varphi$ is normal (in the usual sense) when restricted to any MASA of $\mathcal{A}$.
\end{lem}

Notice that the centre-valued quasitrace on a finite AW*-factor $\A$ is a normal quasitrace. This is because, given any pairwise orthogonal set of projections $(p_i)_{i \in I}$ in $\mathcal{A}$,
    \begin{align*}
         \tau\Bigg(\sup_{\Proj(\mathcal{A})} \{p_i\}\Bigg) =  D\Bigg(\sup_{\Proj(\mathcal{A})} \{p_i\}\Bigg) &= \sum_{i \in I} D(p_i) = \sum_{i \in I} \tau(p_i),
    \end{align*} as required.

\smallskip
We can extend the notion of ``CAP'' or ``normal'' for bounded functions between arbitrary AW*-algebras (or W*-algebras), in the obvious way. 
With this in mind, we state the following auxiliary result which will be used later.

\begin{lem}\cite[Lemma 2]{FeldmanFell} \label{ProjEquivHom}
    Let $\mathcal{A}, \mathcal{B}$ be AW*-algebras and let $\Phi: \mathcal{A} \to \mathcal{B}$ be a $\ast$-homomorphism. If $(p_i)_{i \in I}$ and $(q_i)_{i \in I}$ (indexed by the same set) are two families of pairwise orthogonal projections in $\A$ such that $p_i \sim q_i$ for all $i$, then 
    \begin{align*}
        &\Phi\left(\sup_{\Proj(\mathcal{A})} \{p_i\}\right) = \sup_{\Proj(\mathcal{B})} \Phi(p_i) \\
        \implies &\Phi\left(\sup_{\Proj(\mathcal{A})} \{q_i\} \right) = \sup_{\Proj(\mathcal{B})} \Phi(q_i).
    \end{align*}
\end{lem}

\subsection{Haagerup's Construction}\label{HC}

Continuing the work of Blackadar and Handelman \cite{BH}, Haagerup \cite{Haagerup} introduced a metric on a unital C*-algebra $\mathcal{A}$, which he used to obtain a partial answer to Kaplansky's Conjecture. We introduce this metric now.  

\begin{lem}\cite[Lemma 3.5]{Haagerup}
    Let $\tau$ be a quasitrace on $\mathcal{A}$. Define a map
    \begin{equation*}
        \norm{x}_2 = \tau(x^*x)^{1/2}, \hspace{0.5cm} x \in \mathcal{A}.
    \end{equation*}
    Then
    \begin{enumerate}
        \item $\tau(x+y)^{1/2} \leq \tau(x)^{1/2}+\tau(y)^{1/2}$, for any $x,y \in \mathcal{A}_{+}$.
        \item $\norm{x+y}_2^{2/3} \leq \norm{x}_2^{2/3} + \norm{y}_2^{2/3}$, for any $x, y \in \mathcal{A}$.
        \item $\norm{xy}_2 \leq \norm{x}\norm{y}_2$ and $\norm{xy}_2 \leq \norm{x}_2\norm{y}$ for any $x,y \in \mathcal{A}$.
    \end{enumerate}
\end{lem}

By the preceding lemma, whenever $\tau$ is a \emph{faithful} quasitrace on $\mathcal{A}$, it induces a metric $d_\tau$ on $\A$ given by \[d_\tau(x,y) = \norm{x-y}_2^{2/3}, \hspace{0.5cm} x,y \in \mathcal{A}.\] Notice that if $\tau$ is not faithful, then $d_\tau$ will only be a pseudometric.

Haagerup continues by showing that this metric satisfies a number of useful properties. 

\begin{lem}\cite[Lemma 3.7]{Haagerup}
    Let $\tau$ be a faithful quasitrace on $\mathcal{A}$. Then
    \begin{enumerate} 
        \item The involution $x \to x^*$ is continuous in the $d_\tau$-metric. 
        \item The sum is continuous in the $d_\tau$-metric. 
        \item The product is continuous in the $d_\tau$-metric on bounded subsets of $\mathcal{A}$.
        \item Evaluation $x \to \tau(x)$ is continuous in the $d_{\tau}$-metric on $\mathcal{A}_{+}$.
    \end{enumerate}
\end{lem}

\begin{lem}\cite[Lemma 3.8]{Haagerup} \label{closedBall}
    Let $\tau$ be a faithful quasitrace on $\mathcal{A}$. Then the unit ball of $\A$ is closed in the $d_{\tau}$-metric.
\end{lem}

In the sequel, we will also want to exploit the following result.

\begin{prop}\cite[Proposition 3.10]{Haagerup} \label{completeBall}
    Let $\A$ be a finite AW*-algebra with a normal faithful quasitrace $\tau$. Then the unit ball of $\A$ is complete in the $d_\tau$-metric.
\end{prop}

In order to relate the work of Haagerup and others to the conjectured W*-Pedersen characterization for finite C*-algebras, we want to show that every MASA of a Type II$_1$ AW*-factor is a W*-algebra.
The key idea in the proof (which is really just the proof of \cite[Lemma 3.9]{Haagerup} with a slightly different emphasis) is showing that the unit ball of an AW*-factor $\A$ with unique quasitrace $\tau$ is complete with respect to a particular metric. 
We will need to show that the new metric is appropriately equivalent to $d_\tau$. 
Because the terminology used to describe different metric equivalences is not standard across all references, we clarify the precise definitions below.

\begin{defn}
    Let $d_1, d_2$ be two metrics on the same set $X$. Let ${I: (X,d_1) \to (X, d_2)}$ denote the identity map and let $I^{-1}: (X, d_2) \to (X, d_1)$ be its inverse. 
    \begin{enumerate}[(i)]
        \item We say that $d_1, d_2$ are \emph{equivalent} if both $I$ and $I^{-1}$ are continuous.
        \item We say that $d_1, d_2$ are \emph{uniformly equivalent} if both $I$ and $I^{-1}$ are uniformly continuous.
    \end{enumerate}
\end{defn}

\noindent In general, equivalence of two metrics is not enough to preserve completeness, but uniform equivalence is. 

\begin{prop}\cite[Lemma 3.9]{Haagerup}\label{HaagerupMASA}
    Let $\A$ be a finite AW*-factor equipped with the canonical quasitrace $\tau$. Then each MASA of $\A$ is a W*-algebra.
\end{prop}

\begin{proof}
     Let $M$ be any MASA of $\mathcal{A}$. By Proposition \ref{completeBall}, $\mathcal{A}_1$ is complete with respect to the $d_{\tau}$-metric. 
     By Lemma \ref{closedBall} we know that $\A_1$ is $d_\tau$-closed. 
     In fact, $M_1$ is also $d_\tau$-closed: if $(x_n)_{n \in \mathbb{N}} \subseteq M_1$ is a sequence and $x_n \to x \in A_1$ in the $d_\tau$-metric then for each $m \in M$, continuity of multiplication on bounded sets implies that \[ xm = \lim_{n} x_nm = \lim_{n} m x_n = mx.\]
     Therefore, $x \in M' \cap A = M$, hence $x \in M_1.$
     In particular, $M_1$ must also be complete in the $d_{\tau}$-metric. 
     Since the restriction of the canonical quasitrace $\tau$ to $M$ is a positive, faithful state on $M$, then $d_2(x,y) := \norm{x-y}_2$ is also a metric on $M$. We claim that  $d_\tau$ and $d_2$ are uniformly equivalent metrics on $M_1$. 
     
    We want to show that $I: (M_1, d_\tau) \to (M_1, d_2)$ and $I^{-1}: (M_1, d_2) \to (M_1, d_\tau)$ are both uniformly continuous. Let $\epsilon > 0$ and take $\delta := \epsilon^{2/3}$. Then whenever $x, y \in M_1$ are such that $d_\tau(x,y) < \delta$ we have 
    \begin{align*}
        \epsilon = \delta^{3/2} > d_\tau(x,y)^{3/2} &= \big(\norm{x-y}_2^{2/3}\big)^{3/2} \\
        &= d_2(x,y).
        \end{align*} That is, $I: (M_1, d_\tau) \to (M_1, d_2)$ is uniformly continuous. 
        
    For the other direction, let $\epsilon > 0$ and take $\delta := \epsilon^{3/2}$. Then whenever $x, y \in M_1$ are such that $d_2(x,y) < \delta$ we have  
    \begin{align*}
        \epsilon = \delta^{2/3} > d_2(x,y)^{2/3} &= \norm{x-y}_2^{2/3} \\
        &= d_\tau(x,y),
        \end{align*} so $I^{-1}: (M_1, d_2) \to (M_1, d_\tau)$ is also uniformly continuous. We conclude that $d_{\tau}$ and $d_2$ are uniformly equivalent, hence $M_1$ is complete in the $d_2$-metric.

     By a well-known theorem (see, e.g., \cite[Proposition 2.6.4]{PopaNotes}), completeness of the unit ball of a unital C*-algebra with respect to the $2$-norm induced by a faithful tracial state implies that the  C*-algebra is a W*-algebra. In particular, $M_1$ being $d_2$-complete means that $M$ is a W*-algebra. Since $M$ was an arbitrary MASA of $\mathcal{A}$, we are done.
\end{proof}

\begin{warn} One cannot use Proposition \ref{HaagerupMASA} to show that $\A$ itself is a W*-algebra. While $\mathcal{A}_1$ is complete with respect to $\norm{\cdot}_2$, we cannot apply the cited result because the norm is not \emph{a priori} known to be induced by a faithful tracial state. Indeed, the central problem is that we do not know whether every $2$-quasitrace is a trace. 
\end{warn}

\section{Pedersen-Type Characterizations}\label{PedersenSection}

Before we study the proof of Pedersen's Theorem \ref{PTheorem} and consider the abstract W*-Pedersen Conjecture \ref{absPedersen} in greater detail, let us look at the more general problem of determining the structure of a C*-algebra by the structure of its MASAs.

\subsection{Structure by MASAs}\label{subsec:MASA} In general, demanding that every MASA of a C*-algebra $\A$ satisfies a certain regularity property $(\mathcal{P})$ places very strong requirements on all of $\mathcal{A}$. In fact, demanding that all MASAs of $\A$ satisfy $(\mathcal{P})$ can be \emph{stronger} than requiring that $\A$ satisfies $(\mathcal{P})$.

Recall that a C*-algebra $\A$ is said to have real rank zero (RR0) if arbitrary elements in ${\mathcal{A}}_{sa}$ can be approximated arbitrarily well by elements in ${\mathcal{A}}_{sa}$ with finite spectrum.
In particular, $\A$ has RR0 if for any $a \in \A_{sa}$ and any $\varepsilon > 0,$ there exists some $b \in \A_{sa}$ with finite spectrum such that $\norm{a - b} \leq \varepsilon.$

\begin{exmp}\cite[Theorem 6.1.2]{BlackadarProj}
    Every AW*-algebra has RR0 but a unital C*-algebra with RR0 need not be an AW*-algebra.
\end{exmp}

\begin{prop}
    Let $\A$ be a unital C*-algebra. If every MASA of $\A$ has RR0, then so does $\mathcal{A}$. Conversely, if $\A$ has RR0, its MASAs may not.
\end{prop}

\begin{proof}
    First let $\A$ be a unital C*-algebra in which every MASA has RR0.
    Let $a \in {\mathcal{A}}_{sa}$. 
    We wish to show that $a$ can be approximated arbitrarily well by elements of finite spectrum in ${\mathcal{A}}_{sa}$.
    Let $C^*(1, a)$ be the commutative C*-subalgebra of $\mathcal{A}$ generated by $1$ and $a$. 
    Clearly, $C^*(1, a) \subseteq M$ for some MASA $M$ of ${\mathcal{A}}$.
    Since $M$ has RR0, $a$ can be approximated arbitrarily well by finite spectrum elements in $M_{sa} \subseteq \A_{sa}$.
    Since $a$ was arbitrary, we are done.

    Conversely, let us consider an irrational rotation algebra, 
    \[\A_{\theta} := C^*(u, v: vu = e^{2i\pi\theta}uv),\] where $\theta \in \mathbb{R}\setminus\mathbb{Q}$.
    It is well-known that these algebras have RR0 \cite{ElliottEvans} and that \[M:= C^*(u) \cong C(\mathbb{T})\]is a MASA of $\A_\theta$ (e.g., by \cite{BK26} along with the fact that $\A_{\theta} \cong C(\mathbb{T}) \rtimes_{\theta} \mathbb{Z},$ where the irrational rotation action is \emph{topologically free}).
    Since the real rank of a commutative C*-algebra coincides with the dimension of its underlying topological space \cite[Proposition 1.1]{BrownPedersen}, we see that the real rank of $M$ is $\operatorname{dim}(\mathbb{T}) = 1.$
\end{proof}

\begin{rem}
We do not wish to overstate the extent to which the structure of MASAs of an arbitrary C*-algebra $\A$ determines the full non-commutative structure of $\A$.
As the proposition shows, $\A$ may have a more rigid structure than its MASAs. 

In full generality, it is even possible for $\A$ to have a less rigid structure than its MASAs.
In the 1970s, Akemann and Doner described a non-separable C*-algebra in which every MASA is separable \cite{AkemannDoner}. 
(Those authors constructed this example under the assumption that the Continuum Hypothesis holds; this assumption was later removed by Bice and Koszminder \cite{BiceKoszminder}.)

\end{rem}

Next, we will consider the setting of concrete von Neumann algebras. 
In this case, we are interested in unital C*-algebras whose MASAs are all SOT-closed in the same ambient space $\mathcal{B}(H)$. 

\subsection{Pedersen's Theorem} \label{VNCharacterization}

The following theorem was originally due to Pedersen, \cite{Pedersen}. Because this result does not appear to be well-known, we include the proof as it is presented in \cite{PMonograph}, with some clarifying details added.  

\begin{thm}[Pedersen's Theorem] \label{PTheorem}
    If $\mathcal{A} \subseteq \mathcal{B}(H)$ is a unital C*-algebra such that every MASA is a von Neumann algebra on $H$, then $\A$ is also a von Neumann algebra on $H$. 
\end{thm}

It is clear that any MASA of a von Neumann algebra on $H$ is also a von Neumann algebra on $H$. In this way, it is legitimate (albeit somewhat circular) to characterize von Neumann algebras on $H$ as those C*-subalgebras of $\mathcal{B}(H)$ whose MASAs are von Neumann algebras on $H$. 

\smallskip
We will need some ancillary results in order to  prove Pedersen's Theorem \ref{PTheorem}. 
We only state the first, as the proof follows from a relatively straightforward application of the Stone-Weierstrass theorem along with the usual functional calculus. 

\begin{lem} \cite[Proposition 2.3.2]{PMonograph}\label{strongProp}
    Let $f: \mathbb{R} \to \mathbb{R}$ be a continuous function. If $\abs{f(t)} \leq \alpha \abs{t}+\beta$ for some positive real numbers $\alpha, \beta$, then $f$ is strongly continuous. That is, for a net $(x_i)_{i \in I}$ in $\mathcal{B}(H)_{sa}$ which converges strongly to $x \in \mathcal{B}(H)_{sa}$, the net $\{f(x_i)\}_{i \in I}$ converges strongly to $f(x)$. 
\end{lem}

\begin{rem}
    In the original statement of \cite[Proposition 2.3.2]{PMonograph}, Pedersen additionally assumes that $f(0) = 0$.
    For continuous functions $f: \mathbb{R} \to \mathbb{R}$ satisfying $\abs{f(t)} \leq \alpha \abs{t} + \beta$, but not vanishing at $0$, our more general statement follows by applying that result to the function $f - f(0)$ instead.
\end{rem}

The following technical lemma, also due to Pedersen \cite{PMonograph}, does most of the heavy lifting for the proof of Pedersen's Theorem \ref{PTheorem}.

\begin{lem}\cite[Lemma 2.8.6]{PMonograph} \label{technicalLemma}
    Let $\mathcal{A} \subseteq \mathcal{B}(H)$ be a concrete AW*-algebra. If $p$ is a projection in the SOT-closure of $\mathcal{A}$, then for each pair of vectors $\xi \in pH$ and $\eta \in (1 - p)H$, there is a projection $q_{\xi, \eta}$ in $\A$ such that $q_{\xi, \eta} \xi = \xi$ and $q_{\xi, \eta} \eta = 0$.
\end{lem}

\begin{proof}

    Let $\epsilon > 0$. For each $n \geq 1$, define $\epsilon_n := \frac{1}{2}4^{-n}\epsilon.$ Let $p_0 := 1$ and $x_0 := 0$. 

    By Kaplansky's Density Theorem, there is a net $(u_\lambda)_{\lambda \in \Lambda} \subseteq \mathcal{A}_{1, sa}$ so that $u_\lambda \to p$ strongly. 
   As $(u_\lambda)_{\lambda \in \Lambda}$ is uniformly bounded and converges strongly, we may assume that ${(u_\lambda)_{\lambda \in \Lambda} \subseteq \mathcal{A}_{1,+}}$ by replacing $u_\lambda$ with $u_\lambda^2$, if necessary.
    For our base case, we pick $p_1 := 1$ and choose some $x_1$ from $(u_\lambda)$ such that \[ \norm{(x_1 - p) \xi}^2  = \norm{(x_1 - p)p_1\xi}^2 < \epsilon_1, \quad \norm{(x_1 - p)\xi}^2 < 1, \quad \text{and } \norm{x_1 \eta} < 1.\]

    Suppose that for each $k$, $1 \leq k \leq n$, we have found $x_k \in \{u_\lambda\}$ and $p_k \in \Proj(\mathcal{A})$ such that 
    \begin{enumerate}
        \item $\norm{(x_k - p)p_k\xi}^2 < \epsilon_k, \quad \norm{(x_k - p)\xi} < \frac{1}{k}, \quad$ and $\norm{x_k \eta} < \frac{1}{k}.$
        \item $p_k (x_k - x_{k-1})^2 p_k \leq 2^{-k + 1}.$
        \item $p_k \leq p_{k-1}$ and $\norm{(p_{k-1} - p_k)\xi}^2 \leq 2^{-k+1}\epsilon$.
    \end{enumerate}
    We want to construct $x_{n+1}$ and $p_{n+1}$ satisfying the above properties. 

    Consider the map $\mu_\xi: \mathcal{A} \to \mathbb{C}$ given by $\mu_\xi(a) = \braket{a \xi}{\xi}$ for our fixed $\xi \in H$. 
    By spectral theory, the restriction of $\mu_\xi$ to $C^*(1, p_n (x_n - p)^2 p_n)$ is a bounded Radon measure on the spectrum of $p_n(x_n-p)^2p_n$.
    It follows that we can pick a small open interval $I \subset [2^{-(n+1)}, 2^{-n}]$ and a continuous function $g: \mathbb{R} \to [0,1]$, which satisfies $g(t) = 1$ for all $t \in I$ and 
    \begin{align*}
        \braket{g(p_n (x_n - p)^2 p_n)\xi}{\xi} < \frac{1}{32} \epsilon_n = \frac{1}{8} \epsilon_{n+1}.
    \end{align*}
    
    Now, pick some $t_0 \in I$ and let $\chi_0$ denote the characteristic function on the interval $(-\infty, t_0)$.
    Next, let us choose a function $h: \mathbb{R} \to [0,1]$ such that supp$(h) = \{t \in \mathbb{R}: h(t) \neq 0\} \subseteq I$ and $k = \chi_0+h$ is continuous. 

    For each $\lambda \in \Lambda$, define $y_\lambda:= p_n(u_\lambda-x_n)^2p_n$ and $z_\lambda := p_n(u_\lambda-p)^2p_n$. 
    Then \begin{equation*}
        \chi_0(y_\lambda)z_\lambda \chi_0(y_\lambda) = (k(y_\lambda)-h(y_\lambda))z_\lambda(k(y_\lambda)-h(y_\lambda)).
    \end{equation*} 
    Since $(a+b)^*(a+b) \leq 2a^*a+2b^*b$, we have 
    \begin{equation*} (k(y_\lambda)-h(y_\lambda))z_\lambda(k(y_\lambda)-h(y_\lambda)) \leq 2k(y_\lambda)z_\lambda k(y_\lambda) + 2h(y_\lambda)z_\lambda h(y_\lambda).
    \end{equation*}
    Moreover, $0 \leq {z_\lambda} \leq 4$ and $h^2 \leq g$, so that, combining the equations above, we get 
    \begin{equation*}
        \chi_0(y_\lambda)z_\lambda \chi_0(y_\lambda) \leq 2k(y_\lambda)z_\lambda k(y_\lambda) + 8g(y_\lambda).
    \end{equation*}

    Since $y_\lambda \to p_n(p-x_n)^2p_n$ strongly, $z_\lambda \to 0$ strongly, and the functions $k$, $g$ are continuous, we have that $k(y_\lambda) z_\lambda k(y_\lambda) \to 0$ strongly and $g(y_\lambda) \to g(p_n(p-x_n)^2p_n)$ strongly, by Lemma \ref{strongProp}. 
    In particular, since
    \begin{align*}
        \braket{\chi_0(y_\lambda)z_\lambda \chi_0(y_\lambda) \xi}{\xi} &\leq \braket{2k(y_\lambda)z_\lambda k(y_\lambda)\xi}{\xi} + \braket{8 g(y_\lambda) \xi}{\xi}\\  &\to 8\braket{g(p_n (x_n - p)^2 p_n)\xi}{\xi},
    \end{align*} there exists $\lambda_1 \in \Lambda$ such that \[ \braket{\chi_0(y_\lambda)z_\lambda \chi_0(y_\lambda) \xi}{\xi} < \epsilon_{n+1}\] for all $\lambda \geq \lambda_1$.
    Moreover, since $u_\lambda \to p$ strongly, we can also find $\lambda_2 \in \Lambda$ such that
    \[ \braket{p_n (u_\lambda - x_n)^2 p_n \xi}{\xi} < \epsilon_n, \quad{} \norm{u_\lambda \eta} < \frac{1}{n+1}, \text{ and }  \norm{(u_\lambda - p) \xi} < \frac{1}{n+1}, \] for all $\lambda \geq \lambda_2.$

    Taking $\lambda_3 \in \Lambda$ such that $\lambda_3 \geq \lambda_1$ and $\lambda_3 \geq \lambda_2$, define $p_{n+1} := \chi_0 (y_{\lambda_3}) p_n$ and $x_{n+1} := u_{\lambda_3}$. Clearly $p_{n+1} \leq p_n$, so we have that the pair $x_{n+1}, p_{n+1}$ satisfy property $(1)$.

    Since $\chi_0(t) = 0$ for all $t \geq 2^{-n}$, we have \[ p_{n+1}(x_{n+1} - x_n)^2p_{n+1} = \chi_0(y_{\lambda_3})y_{\lambda_3} \chi_0(y_{\lambda_3}) \leq 2^{-n},\] so property $(2)$ is also satisfied.

    Finally, since $(1-\chi_0) (t) = 0$ for $t < 2^{-(n+1)}$, we have 
    \[ (1 - \chi_0)(y_{\lambda_3}) \leq 2^{n+1} y_{\lambda_3} = 2^{n+1} p_n (x_{n+1} - x_n)^2 p_n,\] which implies
    \[ p_n - p_{n+1} \leq 2^{n+1} p_n (x_{n+1} - x_n)^2 p_n.  \] 
    Using this inequality and the fact that \[\braket{p_n (u_{\lambda_3} - x_n)^2 p_n \xi}{\xi} < \epsilon_n,\] we get
    \[ \braket{(p_n - p_{n+1})\xi}{\xi} \leq 2^{n+1} \epsilon_n = 2^{-n} \epsilon.\] 
    In particular, this means that $p_{n+1}$ satisfies property $(3)$.

    By induction, it follows that we can construct a sequence $(x_n)_{n=1}^{\infty} \in \mathcal{A}_{1,+}$ and a decreasing sequence of projections $(p_n)_{n=1}^{\infty}$ satisfying $(1), (2)$ and $(3)$. 

    Let $\hat{p} := \inf_{\Proj(\mathcal{A})}\{p_n\}$. Since $\A$ is an AW*-algebra, it follows that $\hat{p} \in \mathcal{A}$. Moreover, 
    \[ \braket{(1-\hat{p})\xi}{\xi} = \sum_{n=1}^{\infty} \braket{(p_n - p_{n+1}) \xi}{\xi} \leq \sum_{n=1}^{\infty} 2^{-n} \epsilon = \epsilon, \] where the first equality comes from the fact that we chose $p_1 = 1$.

    By $(2)$, we have that \begin{align*}
        2^{-n} \geq \norm{p_{n+1}(x_{n+1} - x_n)^2 p_{n+1}} 
        &= \norm{((x_{n+1}-x_n)p_{n+1})^*((x_{n+1}-x_n)p_{n+1})} \\
        &= \norm{(x_{n+1} - x_n)p_{n+1}}^2,
    \end{align*} so we know that $(x_n \hat{p})_{n=1}^{\infty}$ converges in norm. Since 
    \begin{align*}
        \norm{x_{n+1}\hat{p}x_{n+1} - x_n \hat{p} x_n} \leq \norm{(x_{n+1}-x_n)\hat{p}x_{n+1}} + \norm{x_n\hat{p}(x_{n+1}-x_n)}  
    \end{align*} and $\norm{x_n} \leq 1$, we conclude that $(x_n \hat{p} x_n)_{n=1}^\infty$ also converges in norm; denote this limit by $y$. 

    By $(1)$, we have $\norm{y \eta} = 0$ and furthermore 
    \begin{align*}
        \braket{y \xi}{\xi} &= \lim \braket{x_n\hat{p}x_n \xi}{\xi} \\
        &= \lim \braket{\hat{p}x_n \xi}{x_n \xi} \\
        &= \braket{\hat{p} \xi}{\xi} \\
        &\geq 1 - \epsilon.
    \end{align*}

    Let $y_{\epsilon}$ denote the support projection of $y$ in $\mathcal{B}(H)$, i.e., the projection onto $\overline{yH}$. 
    Since $y \geq 0$, we have $y_\epsilon = \mathbf{1}_{(0, \norm{y}]}(y)$. Because $\A$ is a concrete AW*-algebra on $H$ it follows that $y_\epsilon$ is a projection in $\mathcal{A}$ such that $y_{\epsilon}\eta = 0$ and $\norm{y_\epsilon\xi}^2 \geq 1 - \epsilon.$

    By repeating the procedure above for each $\epsilon > 0$, we obtain a family of projections 
    ${(y_{\epsilon})_{\epsilon > 0} \subseteq \mathcal{A}}$. 
    Define the projection \[q_{\xi, \eta} := \sup_{\Proj(\mathcal{A})}\{y_{\epsilon}\}.\] Since $\mathcal{A}$ is an AW*-algebra, we have $q_{\xi, \eta} \in \mathcal{A}$. 
    It follows that $q_{\xi, \eta} \xi = \xi$ and $q_{\xi, \eta} \eta = 0$, as desired.
\end{proof}

We are now ready to prove Pedersen's Theorem \ref{PTheorem}.

\begin{proof}[Proof of Pedersen's Theorem.]
    We assume that $\mathcal{A} \subseteq \mathcal{B}(H)$ is a C*-algebra such that every MASA of $\A$ is a von Neumann algebra on $H$. 
    In particular, every MASA of $\A$ is an AW*-algebra and so, by the Sait{\^o}-Wright characterization \ref{AWcharacterization}, $\A$ must be an AW*-algebra.
    
    Let $a \in \mathcal{A}_{sa}$. 
    Since $C^*(1, a)$ is an abelian $C^*$-subalgebra of $\mathcal{A}$, $C^*(1, a)$ is contained in some MASA of $\mathcal{A}$, denoted $M_a$. 
    Since $M_a$ is a von Neumann algebra, it contains each spectral projection of $a$.
    
    Now, notice that any family of orthogonal projections $(p_i)_{i \in I} \subseteq \Proj(\mathcal{A})$ lives in some MASA of $\mathcal{A}$, call it $M$. 
    Since $M$ is a von Neumann algebra, the SOT-sum $\sum_{i \in I} p_i$ is also in $M \subseteq \mathcal{A}$, as required.

    It follows that $\A$ is a concrete AW*-algebra, so we can apply Lemma \ref{technicalLemma}.

    Let $\overline{\mathcal{A}}$ denote the SOT-closure of $\mathcal{A}$. 
    Take any $p \in \Proj ({\overline{\mathcal{A}}})$.  
    Fix $\xi \in pH$ and for each $\eta \in (1-p)H$ construct $q_{\xi, \eta}$ via the procedure in Lemma \ref{technicalLemma}.  
    Let \[q_{\xi} := \inf_{\Proj(\mathcal{A})}\{q_{\xi, \eta}: \eta \in (1-p)H\}.\] 
    Since $\A$ is an AW*-algebra, $q_{\xi} \in \mathcal{A}$. 
    Moreover, it is clear that $q_{\xi}\xi = \xi$ and $q_{\xi} \leq p$. 

    Now let \[q := \sup_{\Proj(\mathcal{A})}\{q_{\xi}: \xi \in pH\}.\]
    Since $q \leq p$ by construction, it suffices to show that we also have $p \leq q$. 
    In fact, this inequality is also immediate from our construction, since $q \xi = \xi$ for all $\xi = p \xi$. 
    In particular, $p=q$, so $p \in \mathcal{A}$. 
    Since $p$ was an arbitrary projection in $\overline{\mathcal{A}}$, we conclude that $\mathcal{A}$ is SOT-closed, hence it is a von Neumann algebra on $H$.
\end{proof}

The following is a well-known result. 
We include it here because the details required to prove it are scattered throughout the literature, but have been given in full above. 
We are not aware of another source in which the necessary ingredients are all gathered together (e.g., it is left as a series of exercises in \cite{Berberian}). 

\begin{cor}\label{AWsub}
    Let $\mathcal{A} \subseteq \mathcal{B}(H)$ be an AW*-algebra. Then $\A$ is a von Neumann algebra on $H$ if and only if $\A$ is an AW*-subalgebra of $\mathcal{B}(H)$.
\end{cor}

\begin{proof}
    The forward direction is obvious. 
    For the converse, recall from Corollary \ref{cor:concreteAW} that $\A$ is an AW*-subalgebra of $\mathcal{B}(H)$ exactly when $\A$ is a concrete AW*-algebra on $\mathcal{B}(H)$.
    Then the argument follows as in the proof of Pedersen's Theorem \ref{PTheorem}.
\end{proof}

We now get the following useful corollary, essentially for free from the preceding result and the relevant definitions.

\begin{cor}\label{cor:CAPrep}
    Let $\mathcal{A}$ be an AW*-algebra. If $(\pi, H)$ is a faithful, CAP representation of $\A$, then $\pi(\A)$ is an AW*-subalgebra of $\mathcal{B}(H)$. In particular, $\A \cong \pi(\A)$ is a W*-algebra. 
\end{cor}

The next consequence of Pedersen's Theorem \ref{PTheorem} is also contained in the original paper \cite{Pedersen}. The details there are sparse, so we follow the proof in \cite{PMonograph}.

\begin{cor}\cite[Theorem 3.9.4]{PMonograph}\label{CAPcor}
    An AW*-algebra $\A$ that admits a separating family of CAP states is a W*-algebra.
\end{cor}

\begin{proof}
    For each CAP state $\omega_i: \mathcal{A} \to \mathbb{C}$, the GNS construction induces a representation $(\pi_i, H_i)$ with a cyclic vector $\xi_i$ such that $\braket{\pi_i(x) \xi_i}{\xi_i} = \omega_i(x)$. 
    Define a pair $(\pi, H) := \oplus_{i} (\pi_i, H_i)$,  where we sum over all the representations of $\mathcal{A}$ coming from the CAP states in the separating family. 
    
    Now take any positive $x \in \ker{\pi}$. 
    For each CAP state $\omega_i$, \[ \omega_i(x) = \braket{\pi_i(x)\xi_i}{\xi_i} = 0,\]
    so we must have that $x = 0$. 
    Therefore, $(\pi, H)$ is a faithful representation of $\mathcal{A}$.

    Now, let $M$ be a MASA of $\mathcal{A}$ and let $\{x_\lambda\}$ be a norm-bounded, monotone increasing net in $M_{sa}$.
    Since $M$ is a commutative AW*-algebra, it is monotone complete. 
    Therefore, $\{x_\lambda\}$ has a LUB, call it $x \in M_{sa}$. 
    As $\{\pi(x_\lambda)\}$ is a monotone increasing net of self-adjoint operators that is bounded above, it converges strongly to some $y \in \mathcal{B}(H)_{sa}$. 
    In fact, $y$ is the LUB for $\{\pi(x_\lambda)\}$ in $\mathcal{B}(H)_{sa}$, so we have that $y \leq \pi(x)$.
    
    For each $\omega_i$ and each unitary $u \in \mathcal{A}$, we have
    \begin{align*}
        \braket{\pi_i(x)\pi_i(u)\xi_i}{\pi_i(u)\xi_i} &= \omega_i(u^*xu) \\
        &= \sup \omega_i(u^* x_\lambda u) \\
        &= \sup \braket{\pi_i(x_\lambda)\pi_i(u)\xi_i}{\pi_i(u)\xi_i} \\
        &= \braket{y\pi_i(u)\xi_i}{\pi_i(u)\xi_i},
    \end{align*}
    where the second equality follows from the fact that a CAP state $\omega_i$ on an AW*-algebra is normal when restricted to any MASA (see Lemma \ref{normalCAP}). 
    
    It follows that $(\pi_i(x)-y)\pi_i(u)\xi_i = 0$ and since $\mathcal{A}$ is spanned by unitaries, we conclude that $(\pi_i(x)-y)[\pi_i(A)\xi_i] = 0$.
    From our construction of $(\pi, H)$, we further conclude that $\pi(x) = y$.

    It follows that the SOT-limit (computed in $\mathcal{B}(H)$) of any monotone increasing net in $\pi(M)_{sa}$ is contained in $\pi(M)_{sa}$. 
    By a result of Kadison \cite[Lemma 1]{Kadison}, this is enough to conclude that $\pi(M)$ is a von Neumann algebra on $H$.
    Since $M$ was arbitrary, it follows that each MASA of $\pi(\mathcal{A})$ is a von Neumann algebra on $H$. 
    By Pedersen's Theorem \ref{PTheorem}, we conclude that $\pi(\mathcal{A})$ is a von Neumann algebra on $H$. 
    In other words, $\A$ is a W*-algebra.
\end{proof}

\begin{rems} The preceding corollary could have been used in several places above.
\begin{enumerate}
    \item  Instead of the argument adopted in Corollary \ref{AWsub}, we could use Corollary \ref{CAPcor} to prove that if $\A$ is an AW*-subalgebra of $\mathcal{B}(H)$ then $\A$ is a von Neumann algebra on $H$. 
    This follows by noticing that, in this case, \[\{\braket{(\cdot)\xi}{\xi}: \xi \in H,\, \norm{\xi} = 1\}\] is a separating family of CAP states for $\mathcal{A}$. 
    \item  Instead of arguing as Haagerup did  in Proposition \ref{HaagerupMASA}, to see that every MASA of a Type II$_1$ factor is a W*-algebra, we could have proceeded as follows.
    Note that the restriction of the quasitrace $\tau$ to each MASA $M$ of a Type II$_1$ AW*-factor is a faithful state. 
    Then applying Corollary \ref{CAPcor} to the separating (singleton) family $\{\tau|_M\}$ it is immediate that each MASA $M$ is a W*-algebra. 
\end{enumerate}
    
\end{rems}

\subsection{W*-Pedersen Conjecture}\label{WPedSection}

We now move towards a less concrete version of Pedersen's Theorem, beginning with the following corollaries. 
While these first few results may be known to some experts, to the best of our knowledge, none have been explicitly stated or proved in the literature.  

\smallskip
In this section, we will consider the notion of a \emph{W*-subalgebra} of a W*-algebra $\B$. By this, we mean a unital weak*-closed C*-subalgebra of $\B$ whose W*-structure agrees with that of $\B$.

\begin{cor}\label{CorW}
    Let $\A$ be a unital C*-algebra and let $\B$ be a W*-algebra. If $\mathcal{A} \subseteq \mathcal{B}$ is a C*-subalgebra and every MASA of $\mathcal{A}$ is a W*-subalgebra of $\mathcal{B}$, then $\A$ is a W*-algebra. 
\end{cor}

\begin{proof}
    As above, the fact that every MASA of $\A$ is a W*-algebra implies that $\A$ is an AW*-algebra.
    Let $(\pi, H)$ be a normal, faithful representation of $\mathcal{B}$, so that $\pi(\mathcal{B}) \subseteq \mathcal{B}(H)$ is a von Neumann algebra. 
    Since every MASA of $\pi(\mathcal{A})$ is of the form $\pi(M)$ where $M$ is a MASA of $\mathcal{A}$, it suffices to show that each $\pi(M)$ is a von Neumann algebra on $H$.
    As normal $\ast$-homomorphisms preserve W*-algebraic structure, we know that $\pi(M)$ is a W*-subalgebra of $\pi(B)$ and therefore a von Neumann algebra on $H$. 
    So $\pi(\mathcal{A}) \subseteq \mathcal{B}(H)$ is a unital C*-algebra in which every MASA is a von Neumann algebra on $H$ and the result follows from Pedersen's Theorem \ref{PTheorem}.
\end{proof}

\begin{cor}\label{cor:AWsubalgW}
    Let $\A$ be an AW*-algebra and let $\B$ be a W*-algebra. If $\A$ is an AW*-subalgebra of $\B$ then $\A$ is a W*-algebra. 
\end{cor}

\begin{proof}
    As above, let $(\pi, H)$ be a normal, faithful representation of $\mathcal{B}$, so that $\pi(\mathcal{B}) \subseteq \mathcal{B}(H)$ is a von Neumann algebra. 
    Since normal $\ast$-homomorphisms preserve W*-algebraic structure, they must also preserve AW*-algebraic structure.
    It follows that $\pi(\mathcal{A})$ is an AW*-subalgebra of $\pi(\mathcal{B})$ and moreover, $\pi(\mathcal{A})$ is an AW*-subalgebra of $\mathcal{B}(H)$.
    By Corollary \ref{AWsub}, $\pi(\mathcal{A})$ is a von Neumann algebra on $H$, so $\mathcal{A}$ is a W*-algebra.
    In particular, $\A$ is a W*-subalgebra of $\mathcal{B}$.
\end{proof}

Putting the two preceding corollaries together, we get the following.

\begin{cor}\label{cor:abstractAnalogue}
    Let $\A$ be an AW*-algebra, let $\B$ be a W*-algebra, and let $\mathcal{A} \subseteq \mathcal{B}$ be an inclusion of unital C*-algebras. If every MASA of $\A$ is an AW*-subalgebra of $\mathcal{B}$, then $\A$ is a W*-algebra. In particular, $\A$ is a W*-subalgebra of $\mathcal{B}$.
\end{cor}

Notice that, while Corollary \ref{cor:abstractAnalogue} is the exact W*-analogue of Pedersen's Theorem \ref{PTheorem}, it is not fully intrinsic -- it relies on the realization of the AW*-algebra as an AW*-subalgebra of a W*-algebra. 
We conjecture that the following intrinsic characterization holds.

\begin{conj}[W*-Pedersen]\label{absPedersen}
    A unital C*-algebra $\A$ is a W*-algebra if and only if every MASA of $\A$ is a W*-algebra. 
\end{conj}

We now offer some evidence that the W*-Pedersen characterization is true, so that it is a legitimate conjecture (rather than simply an unjustified statement). 
First, it is clear that if $\A$ is a W*-algebra, then every MASA of $\A$ is also a W*-algebra, so the forward direction holds.

Partial evidence for the converse follows from the Sait{\^o}-Wright characterization \ref{AWcharacterization} of AW*-algebras: any unital C*-algebra $\A$ with W*-MASAs is an AW*-algebra with some extra W*-structure on its commutative parts. 
So, the open question becomes the following.

\begin{quest}\label{q:commW}
    Is the extra structure that turns an AW*-algebra into a W*-algebra determined commutatively?
\end{quest}

We have hope that this question has a positive resolution from the considerations of Section \ref{subsec:MASA}, which show that placing a strong regularity property on the MASAs of a C*-algebra is often very restrictive. 
Moreover, it is obvious that the commutative structure is enough to determine the W*-structure in the case of a commutative C*-algebra. 
In fact, it is already known that this is (more than) enough for any Type I AW*-algebra.

\begin{prop}\label{TypeI}
    Let $\mathcal{A}$ be a Type I AW*-algebra. Then $\A$ is a W*-algebra if and only if every MASA is a W*-algebra.
\end{prop}

\begin{proof}
    We only need to show the converse.
    If each MASA $M_i$ of $\A$ is a W*-algebra then $\mathcal{Z} (\mathcal{A}) = \cap_{i \in I} M_i$ is also a W*-algebra. Note that this follows because $\mathcal{Z}(\mathcal{A})$ is an AW*-subalgebra of each of the W*-algebras $M_i$ so by Corollary \ref{cor:AWsubalgW}, $\mathcal{Z}(\A)$ is a W*-algebra.
    
    Moreover, by \cite[Theorem 2]{Kap52}, an AW*-algebra of Type I with W*-algebra centre is itself a W*-algebra.
\end{proof}

\begin{rem}
    The key idea in our proof of Proposition \ref{TypeI} -- which will be a recurring theme throughout -- is the fact that if every MASA of $\A$ is a W*-algebra, then $\mathcal{Z(A)}$ is a W*-algebra. 
    But the latter is a strictly weaker condition: there are AW*-algebras with W*-centre whose MASAs are not all W*-algebras (e.g., the Type III non-W*, AW*-factors constructed in \cite{Takenouchi} and \cite{Dyer}).  
\end{rem}

In light of Corollary \ref{CAPcor}, one method of attack for showing that the W*-Pedersen characterization holds in general, would be to consider a separating collection of normal states on each MASA of $\A$ and try to extend them to (enough) normal states on all of $\mathcal{A}$. 
Perhaps unsurprisingly, this seems a difficult thing to do. 
It requires the use of some Mackey-Gleason type results, which are known for von Neumann algebras, but are not known for general AW*-algebras (see Chapter 5, \cite{HamTextbook} for the Generalized Gleason Theorem; see \cite{Hamhalter}, by that same author, for some extensions to the AW*-algebra case). 
When we restrict to the properly infinite case, we can use these (and other) known results to make progress towards new characterizations of properly infinite W*-algebras. We do exactly this in Section \ref{ProperlyInfiniteSec}.

Given the lack of Mackey-Gleason results for arbitrary finite AW*-algebras, we cannot prove that the W*-Pedersen characterization holds for finite C*-algebras.
However, we \emph{can} prove that showing the W*-Pedersen characterization holds for finite C*-algebras is equivalent to showing that every $2$-quasitrace is a trace. 

\section{Equivalences} \label{MainSection}

\subsection{Finite W*-Pedersen implies Kaplansky}

We are now ready to show that if we assume the W*-Pedersen characterization holds for finite C*-algebras, then any Type II$_1$ AW*-factor is a W*-algebra. 

\begin{prop} \label{WPtoKaplansky} Suppose that the W*-Pedersen characterization holds for finite C*-algebras. If $\A$ is a Type II$_1$ AW*-factor, then it is a W*-algebra. 
\end{prop}

\begin{proof}
    Any Type II$_1$ AW*-factor is, by definition, a finite AW*-factor. By Proposition \ref{HaagerupMASA}, every MASA of $\mathcal{A}$ is a W*-algebra. 
   Under the assumption that the W*-Pedersen characterization \ref{absPedersen} holds for finite C*-algebras, then $\A$ is a W*-algebra. 
   In particular, it is a Type II$_1$ W*-factor and its unique quasitrace is a bona fide trace.
\end{proof}

\subsection{Kaplansky implies 2Q} 

Through deep results in the works of Blackadar and Handelman \cite{BH} and Haagerup \cite{Haagerup}, it is already known that a positive solution to Kaplansky's Conjecture \ref{KapConjecture} exists if and only if a positive solution to the 2Q Conjecture \ref{2Q}. In this subsection, we endeavour to convince the reader of the forward direction of this equivalence by sketching the relevant proof(s).

Recall that $QT(\mathcal{A})$ is a convex set that is compact in the topology of pointwise convergence. We will call an extreme point of $QT(\mathcal{A})$ an \emph{extremal quasitrace}. For any $\tau \in QT(\mathcal{A})$, we can define the \emph{trace-kernel ideal} of $\tau$ in $\mathcal{A}$ to be the set \[ I_{\tau} := \{a \in \mathcal{A}: \tau(a^*a) = 0\}. \] 
 It is clear that if $\tau$ is faithful on $\mathcal{A}$, then $I_{\tau} = \{0\}$. 
On the other hand,  if $\tau$ is not faithful on $\mathcal{A}$, we would like to pass to the quotient algebra $\mathcal{A}/I_{\tau}$. 
It turns out that, even though $\tau$ is not \emph{a priori} linear, $I_{\tau}$ is indeed a closed (two-sided) ideal in $\mathcal{A}$. 
Moreover, there is a unique faithful quasitrace $\bar{\tau}$ on $\mathcal{A}/I_{\tau}$ such that $\tau(x) = (\bar{\tau} \circ q) (x)$, where $q$ denotes the quotient map (see Proposition 3.2, \cite{Haagerup}).

Finally, in order to pass from arbitrary unital C*-algebras with quasitraces to AW*-algebras, we will need another result and a useful definition/construction.

\begin{prop}\cite[Corollary 4.3]{Haagerup}\label{morphism}
   Let $\A$ be a unital C*-algebra with a faithful quasitrace $\tau$. Then there is an injective $\ast$-homomorphism $\pi$ of $\A$ into a finite AW*-algebra $\mathcal{M}$ with a faithful normal quasitrace $\overline{\tau}$ such that
   \begin{equation*}
       \tau(a) = \overline{\tau} \circ \pi(a), \hspace{0.5cm} a \in \mathcal{A}. \end{equation*}
\end{prop}

\begin{defn}
     Let $(\mathcal{A}, \tau)$ and $(\mathcal{M}, \overline{\tau})$ be as in Proposition \ref{morphism}. We call the $d_{\bar{\tau}}$-closure of $\pi(\A)$ in $\mathcal{M}$ the \emph{AW*-completion} of $\mathcal{A}$.\footnote{This definition of the AW*-completion differs from the original \cite[Definition 4.5]{Haagerup}. We find this description of the AW*-completion more manageable for our purposes. The AW*-algebra we obtain from the definition is $\ast$-isomorphic to the one coming from the original (see \cite[Proposition 3.10]{MThesis}).}
\end{defn}

As the name suggests, the AW*-completion of a unital C*-algebra $\A$ with a faithful quasitrace $\tau$ is always an AW*-algebra (see \cite[Section 4]{Haagerup}). Moreover, if $\tau$ is an extreme point of $QT(\mathcal{A})$ then the AW*-completion of $\A$ is a finite AW*-factor (see \cite[Proposition 4.6]{Haagerup}).

\begin{prop}\label{Kaplanskyto2Q}
    If every Type II$_1$ AW*-factor is a W*-factor, then every $2$-quasitrace on a unital C*-algebra is a trace. 
\end{prop}

\begin{proof}[Sketch of proof.]
    Let $\mathcal{A}$ be a unital C*-algebra which admits a $2$-quasitrace $\tau$.
    First, let us suppose that $\tau$ is extremal.
    By taking the quotient $\mathcal{A}/I_{\tau}$ if necessary, we may also assume that $\tau$ is faithful.
    The AW*-completion of $(\A, \tau)$ is a finite AW*-factor with its unique centre-valued quasitrace $(\mathcal{M}, \overline{\tau})$.

    If $\mathcal{M}$ is Type I, it is automatically a W*-algebra by \cite[Theorem 2]{Kap52}.
    If it is a Type II$_1$ AW*-factor then by our assumption, it is also a W*-factor. 
    In either case, $\mathcal{M}$ can be faithfully represented as a von Neumann algebra, so $\overline{\tau}$ is a trace.
    Since the restriction of $\overline{\tau}$ to $\mathcal{A}$ is exactly $\tau$, the latter must also be a trace.
    If $\tau$ is not an extremal quasitrace then, by the Krein-Milman Theorem, $\tau$ must be arbitrarily close to a convex combination of extremal quasitraces (i.e., traces). 
    By a standard argument, such a $\tau$ must also be linear.
\end{proof}

In the statement above, it is important to emphasize that we are considering $2$-quasitraces. 
In 2006, in an unpublished manuscript, Kirchberg gave an example of a $1$-quasitrace on a unital C*-algebra that is not a trace. This example can be found in \cite{Kirchberg}, an unfinished manuscript written by Kirchberg, which was made available posthumously.

\subsection{2Q implies Finite W*-Pedersen}

If we assume that every $2$-quasitrace is a trace, then clearly every Type II$_1$ AW*-factor has a trace (and therefore each Type II$_1$ AW*-factor is a Type II$_1$ W*-factor). 
However,  we would like to show that a positive solution to the 2Q Conjecture implies that the W*-Pedersen characterization holds for \emph{any} finite C*-algebra, not just those with trivial centre.
For this, we first want to move from the setting of Type II$_1$ AW*-factors to arbitrary Type II$_1$ AW*-algebras. 

\begin{lem}\label{lem:linearQT}
    Suppose that any $2$-quasitrace on a unital C*-algebra is a trace. Then the unique centre-valued quasitrace on a Type II$_1$ AW*-algebra is linear. 
\end{lem}

\begin{proof}
    Let $\A$ be a Type II$_1$ AW*-algebra with unique centre-valued quasitrace  $T$. 
    Let $\tau$ be a genuine quasitrace on $\mathcal{A}$ (that is, $\tau: \mathcal{A}\to\mathbb{C}$ is a normalized $2$-quasitrace).
    Recall that any quasitrace on $\A$ can be uniquely written as $\tau = \varphi \circ T$, for some state $\varphi:\mathcal{Z(A)} \to \mathbb{C}$ (see e.g., Theorem II.1.7, \cite{BH}).
    By assumption, $\tau = \varphi \circ T$ is linear. 
    Let us pick arbitrary $x,y \in \mathcal{A}_{+}$. Then
    \begin{align*}
        \tau(x+y) &= \tau(x)+\tau(y) \\
        \iff \varphi(T(x+y)) &= \varphi( T(x) + T(y)).
    \end{align*}
    Let $z := T(x+y)-T(x)-T(y) \in \mathcal{Z(A)}$. It follows from the computation above that $\varphi(z) = 0$.
    Moreover, it is clear that for any state $\psi:\mathcal{Z}(A) \to \mathbb{C}$, the map $\psi \circ T$ is also a genuine quasitrace on $\mathcal{A}$. 
    As argued above, we can therefore deduce that $\psi(z) = 0$ for any state $\psi:\mathcal{Z}(A) \to \mathbb{C}$.
    In particular, this means that $z = 0$.
    Since $x,y \in \mathcal{A}_{+}$ were arbitrary, this means that the centre-valued quasitrace $T$ is linear on $\mathcal{A}_{+}$.
    By a standard argument, it follows that $T$ is linear on all of $\mathcal{A}$.
\end{proof}

Now we are ready to show the next implication.

\begin{prop}\label{2QtoWP}
    Suppose that any quasitrace on a unital C*-algebra is a trace. Then the W*-Pedersen characterization holds for all finite C*-algebras.
\end{prop}

\begin{proof}
    Let $\A$ be a finite C*-algebra in which every MASA is a W*-algebra.
    We denote the set of all MASAs by $\{M_i: i \in I\}.$
    By the Sait{\^o}-Wright characterization \ref{AWcharacterization}, $\A$ is an AW*-algebra.
    Since we can decompose $\A$ into Type I and Type II$_1$ parts, i.e., $\mathcal{A} = \mathcal{A}_I \oplus \mathcal{A}_{II_1}$, we break into separate cases.
    
    First, if we assume that $\A$ is Type I, then the fact that $\A$ is a W*-algebra is immediate from the fact that $\mathcal{Z(A)} = \cap_{i \in I} M_i$ is a W*-algebra (see Proposition \ref{TypeI}).
    
    Now, let us assume that $\A$ is Type II$_1$. 
    Let $T$ be the canonical centre-valued quasitrace on $\mathcal{A}$; under our assumptions and Lemma \ref{lem:linearQT}, $T$ is linear.
    Then $T$ is a faithful CAP conditional expectation onto the W*-algebra $\mathcal{Z(A)}$ (notice that it is a projection of norm one).
    
    Pick any $0 \neq a \in \mathcal{A}_{+}$.
    Since $T$ is positive, $T(a) \in \mathcal{Z(A)}_{+}$.
    Moreover, since $\mathcal{Z(A)}$ is a W*-algebra, there is some CAP state $\varphi_a \in \mathcal{Z(A)}_*$ such that $\varphi_a(T(a)) \neq 0$.
    Then ${(\varphi_a \circ T)_{a \in \mathcal{A}_{+} \setminus\{0\}}}$ defines a separating family of CAP states.
    So by Corollary \ref{CAPcor}, $\A$ is a W*-algebra.

    The direct sum of W*-algebras is a W*-algebra, so $\A = \A_I \oplus \A_{II_1}$ is a W*-algebra.
\end{proof}

\subsection{Retracts of Biduals}\label{bidualChar} In this section, we establish the final equivalence of our main theorem. In particular, we show that the final condition is directly equivalent to a positive solution to Kaplansky's Conjecture. In order to do this, we will need the following result.

\begin{lem}\label{quotientLinear}
    Let $\A$ be a Type II$_1$ AW*-algebra with $\mathcal{Z(A)}$ a W*-algebra. Then $\A$ is a W*-algebra if and only if there exists a W*-algebra $\mathcal{M}$, without Type I$_2$ direct summand, and a norm-closed two-sided ideal $I$ of $\mathcal{M}$ such that $\mathcal{A} \cong \mathcal{M}/I$.
\end{lem}

\begin{proof}
    If $\A$ is a Type II$_1$ W*-algebra then $\mathcal{A} \cong \mathcal{A}/\{0\}$, as required.
    
    On the other hand, suppose that we have a W*-algebra $\mathcal{M}$ without Type I$_2$ direct summand and a norm-closed two-sided ideal $I$ of $\mathcal{M}$ such that $\mathcal{A} = \mathcal{M}/I$.
    Recall that the Mackey-Gleason Theorem says that any positive, bounded quasilinear functional on $\mathcal{M}$ is linear; moreover, Bunce and Wright showed (see \cite[Corollary]{BunceWright}) that this result extends to quotients, so any quasilinear functional on $\mathcal{A} \cong \mathcal{M}/I$ must be linear.
    
    If $\A$ is a factor, then these results imply that the canonical quasitrace on $\A$ is linear. In particular, since the canonical quasitrace is CAP and faithful, $\A$ must be a Type II$_1$ W*-factor (cf. Corollary \ref{CAPcor}).
    Otherwise, for the non-factorial case, let $T$ denote the centre-valued quasitrace on $\mathcal{A}$ and recall that $\mathcal{Z(A)}$ admits a separating family of CAP states $(\varphi_j)_{j \in J} \subseteq \mathcal{Z}(\A)_*$, as it is a W*-algebra. 
    It is not hard to see that each $\varphi_j \circ T$ is a CAP quasistate on $\mathcal{A}$, hence it is linear by the cited results. 
    In particular, $\A$ admits a separating family of CAP states, making it a W*-algebra, by Corollary \ref{CAPcor}.
\end{proof}

In general, a unital C*-algebra that can be realized as the C*-quotient of a W*-algebra does not need to be a W*-algebra (e.g., the Calkin algebra). 
In fact, there are even AW*-algebras that can be realized in this way, but are not W*-algebras.
For the sake of exposition, we relegate further consideration of C*-quotients of W*-algebras to Section \ref{MCQuestion}.
That discussion provides additional context and motivation for the following definition.

\begin{defn}\label{defn:retract}
    Let $\A$ be a unital C*-algebra. We say that $\A$ is a \emph{($\ast$-homomorphic) retract of its bidual} if there exists a multiplicative conditional expectation $E: \mathcal{A}^{**} \twoheadrightarrow \mathcal{A}$ (that is, a $\ast$-homomorphism which fixes elements of $\mathcal{A}$). In this case, $\mathcal{A} \cong \mathcal{A}^{**}/\ker{E}$.
\end{defn}

\begin{exmp}\label{bidualRetractEx}
    We show that any W*-algebra $\mathcal{M}$ is a retract of its bidual. 
    
    Let $\mathcal{M}_*$ be the predual of $\mathcal{M}$ and $\iota: \mathcal{M}_* \to (\mathcal{M}_*)^{**}$ the canonical embedding of the predual into its bidual. 
    We can view the dual map $\iota^{*}: (\mathcal{M}_*)^{***} \to (\mathcal{M}_*)^{*}$ as a map $E: \mathcal{M}^{**} \twoheadrightarrow \mathcal{M}$ under the usual identifications.
    In particular, there is a central projection $z \in M^{**}$ such that $E$ is the canonical projection of $M^{**}$ onto $zM^{**} \cong M$, which is a multiplicative conditional expectation.
\end{exmp}

\begin{cor}\label{cor:bidualII_1}
    Let $\A$ be a Type II$_1$ AW*-factor. Then $\A$ is a W*-factor if and only if it is a retract of its bidual.
\end{cor}

\begin{proof}
    If $\A$ is a W*-algebra, then by Example \ref{bidualRetractEx}, it is a retract of its bidual.
    
    To prove the converse, let us suppose that $\A$ is a retract of its bidual. 
    In other words, we suppose that $\mathcal{A} \cong \A^{**}/\ker{E}$. 
    Notice that by Lemma \ref{quotientLinear}, it suffices to show that if $\A$ is a Type II$_1$ AW*-factor then $\mathcal{A}^{**}$ does not have any Type $I_2$ direct summand. 
    
    Indeed, if $\A^{**}$ had a Type I$_2$ direct summand, then there would be some central projection $z \in \Proj(\A^{**})$ such that \[z\A^{**} \cong M_2(\mathcal{Z}(z\A^{**})).\] 
    Let $\chi: \mathcal{Z}(z\A^{**}) \to \mathbb{C}$ be any character.
    Composing the map $a \mapsto za$ with the isomorphism above and the homomorphism on $M_2(\mathcal{Z}(z\A^{**}))$ induced entrywise by $\chi$ gives a unital, finite-dimensional representation of $\A$. 
    By simplicity of $\mathcal{A}$, this representation would be faithful, forcing $\A$ to be finite-dimensional, which is impossible.    
\end{proof}

Putting Propositions \ref{WPtoKaplansky}, \ref{Kaplanskyto2Q}, and \ref{2QtoWP} together with Corollary \ref{cor:bidualII_1}, we arrive at our first main theorem.

\begin{thm}\label{MainThm}
    The following are equivalent:
    \begin{enumerate}[(i)]
        \item Every Type II$_1$ AW*-factor is a Type II$_1$ W*-factor.
        \item Any $2$-quasitrace on a unital C*-algebra is a trace.
        \item A finite C*-algebra is a W*-algebra if and only if all of its MASAs are W*-algebras.
        \item Every Type II$_1$ AW*-factor is a retract of its bidual.
    \end{enumerate}
\end{thm}

\begin{rem}\label{finiteRem}
    For the implications involving the W*-Pedersen characterization for finite C*-algebras $\mathcal{A}$, notice that we essentially only made use of the fact that $\mathcal{Z(A)}$ is a W*-algebra. 
Indeed, we could replace $(iii)$ above by 
\begin{enumerate}[(iii')]
    \item A finite C*-algebra is a W*-algebra if and only if all of its MASAs are AW*-algebras and its centre is a W*-algebra.
\end{enumerate}
However, in the non-finite case, condition (iii') is not enough to characterize W*-algebras. 
Indeed, it is known that there are even Type III AW*-factors that are not W*-factors (see e.g., \cite{Dyer}, \cite{Takenouchi}).
Ultimately, we conjecture that the W*-Pedersen characterization holds in the properly infinite case, as well,  so we leave the condition as stated. 
\end{rem}

\subsection{Spatial Characterization of II$_1$ Factors}

The statement of Theorem \ref{MainThm} highlights the difficulty in resolving Kaplansky's Conjecture or the 2Q Conjecture. 
It is equivalent to establishing a global property (that a finite C*-algebra is a W*-algebra) from only piecewise local properties (that each MASA is a W*-algebra). 
If we have access to topological structure on the MASAs, as in the concrete Pedersen Theorem \ref{PTheorem}, the local-to-global problem can be handled; without the additional topological structure, it is difficult.

Nevertheless, we are able to establish a result which \emph{does} characterize the Type II$_1$ AW*-factors that are W*-factors and involves a spatial condition which is \emph{a priori} weaker than the one in Pedersen's Theorem \ref{PTheorem}. 
In order to proceed, we need the following technical result. 

\begin{lem}\label{diffuseMASA}
    Let $\A$ be a Type II$_1$ AW*-factor with its canonical quasitrace $\tau$ and let $M$ be any MASA of $\mathcal{A}$. For every projection $q \in M$ and every $t \in [0, \tau(q)]$, there exists a projection $q' \in M$ such that $q' \leq q$ and $\tau(q') = t$.
\end{lem}

\begin{proof}
    We first show that every MASA $M$ of $\mathcal{A}$ is diffuse (i.e., there are no minimal projections in $M$). 
    Suppose, by way of contradiction, that $p \in M$ is a minimal projection. Then $pMp = \mathbb{C}p$. 
    Since $\A$ is Type II$_1$, it contains no abelian projections. 
    In particular, $p\mathcal{A}p \neq \mathbb{C}p$ and so there exists a projection $q \in p\mathcal{A}p$ with $0 \neq q \lneq p$.
    But then for any $m \in M$, since $mp = \lambda_mp$ for some scalar $\lambda_m$, we get \[ qm = qpm = q(mp) = q\lambda_mp = \lambda_mqp = \lambda_mpq = mpq = mq.\] 
    Then $C^{*}(M, q)$ is a commutative C*-subalgebra of $\A$ containing $M$, which contradicts that $M$ is a MASA.

    For every projection $q \in M$ and every real number $t \in [0, \tau(q)]$, we can use the fact that $M$ is diffuse to prove that we can find another projection $q' \in M$ such that $q' \leq q$ and $\tau(q') = t$.

    Let $\mathcal{P}_t := \{p \in \Proj(M): p \leq q, \,\tau(p) \leq t\}$. 
    By the fact that $\tau$ is normal on $M$, we can use a standard Zorn's Lemma argument\footnote{If $\mathcal{C} \subseteq \mathcal{P}_t$ is a chain, then there is some supremum $p_{\mathcal{C}} = \sup_{\Proj(M)} \{r: r \in \mathcal{C}\}$ and $\tau(p_{\mathcal{C}}) = \sup_{r \in \mathcal{C}} \tau(r)\leq t$.} to see that $\mathcal{P}_t$ admits a maximal element, call it $p$, with $\tau(p) \leq t.$

    If $\tau(p) < t \leq \tau(q),$ then $e:= q-p$ is a non-zero projection in $M$.
    Since $M$ is diffuse, then $e$ is not minimal, so there is another non-zero proper projection $e_1 < e.$
    Writing $e = e_1 + (e-e_1),$ we see that we must have $$\tau(r_1) \leq \frac{\tau(e)}{2}$$ where $r_1 := e_1$ or $r_1 := e - e_1.$ We fix $r_1$ so that the inequality holds.

    Repeating this argument, we find non-zero projections $r_k \leq r_{k-1} < e$ with $$\tau(r_k) \leq \frac{\tau(e)}{2^k}.$$ 
    After enough divisions of $\tau(e)$, we get a projection $0 \neq r_n \leq q - p$ such that $$\tau(r_n) \leq \frac{\tau(e)}{2^n} \leq t- \tau(p).$$
    Therefore, $p+r_n$ is a projection with $p + r_n \leq q$ and $$\tau(p+r_n) = \tau(p) + \tau(r_n) \leq t,$$
    so $p+r_n \in \mathcal{P}_t.$ 
    
    At the same time, we have $p + r_n > p$, contradicting the maximality of $p \in \mathcal{P}_t.$
    Therefore, we must have $\tau(p) = t.$
    Taking $q' := p$ gives the desired result.   
\end{proof}

The following argument was pointed out to the author by Hannes Thiel, who noticed that the author's argument in an earlier version of Proposition \ref{TypeIIIFactor} could be extended. 

\begin{thm}\label{minimalPedersenII1}
    Let $\A$ be a Type II$_1$ AW*-factor. Then $\A$ is a W*-factor if and only if there exists a $\ast$-representation $(\pi, H)$ of $\A$ and a MASA $M$ of $\mathcal{A}$ such that $\pi|_M$ is non-zero and CAP.
\end{thm}

\begin{proof}
    We need only prove the converse. 
    First, since $\A$ is a simple C*-algebra, either $\pi$ is the zero $\ast$-homomorphism or it is faithful.
    If $\pi$ were zero, then we would have $\pi(m) = 0$ for all $m \in M$, contradicting our assumption. 
    Therefore, $\pi$ is faithful and so $\mathcal{A} \cong \pi(\mathcal{A})$.
    
    Moreover, since $\pi|_M$ must also be faithful, we have that $\pi|_M: M \to \mathcal{B}(H)$ is a faithful, CAP $\ast$-homomorphism between W*-algebras, hence ultraweakly continuous, and so $\pi(M)$ is a von Neumann algebra on $H$. 
    We wish to show that for \emph{any} MASA $N$ of $\mathcal{A}$, $\pi(N)$ is also a von Neumann algebra on $H$, so that we can invoke the concrete Pedersen Theorem \ref{PTheorem}. 

    Fix a MASA $N \subseteq \mathcal{A}$. Let $(p_i)_{i \in I}$ be a family of pairwise orthogonal projections in $N$. 
    If $I$ is finite, then clearly \[\pi\left(\sum_{i} p_i \right) = \sum_i \pi(p_i).\] 
    Since $\tau|_N$ is faithful and normal, $N$ is countably decomposable and therefore, without loss of generality, we may assume that $I = \mathbb{N}$. 
    Thus, let $(p_k)_{k \in \mathbb{N}}$ be a sequence of pairwise orthogonal projections in $N$. 
    We will inductively choose a sequence of pairwise orthogonal projections $q_k \in M$ (the MASA on which we know $\pi$ to be completely additive) such that $q_k \sim p_k$ for each $k$.

    By Lemma \ref{diffuseMASA} with $t := \tau(p_1) \leq 1 = \tau(1)$, we can find a projection $q_1 \in M$ with $\tau(q_1) = \tau(p_1)$. Since equivalence of projections is determined by $\tau$, we have $q_1 \sim p_1$. 

    Next, since $\tau$ is additive on $M$ and on $N$, we have \[ \tau(p_2) \leq 1 - \tau(p_1) = 1 - \tau(q_1) = \tau(1-q_1).\] 
    Applying Lemma \ref{diffuseMASA} again, we obtain a projection $q_2 \in M$ with $q_2 \leq (1-q_1)$ and $\tau(q_2) = \tau(p_2).$ In other words, $q_1$ and $q_2$ are orthogonal and $q_2 \sim p_2$.

    Assume we have constructed pairwise orthogonal projections $q_1, \ldots, q_n \in M$ with $p_j \sim q_j$ for $j = 1, \ldots, n$.
    Then \[ \tau(p_{n+1}) \leq 1 - \sum_{k=1}^{n} \tau(p_k) = 1 - \sum_{k=1}^n \tau(q_k) = \tau\left(1 - \left(\sum_{k=1}^{n} q_k\right)\right).\]
    Again, by applying Lemma \ref{diffuseMASA}, we obtain a projection $q_{n+1} \in M$ with \[q_{n+1} \leq \left(1-\sum_{k=1}^n q_k\right)\] and \[\tau(q_{n+1}) = \tau(p_{n+1}).\]
    In other words, we have pairwise orthogonal projections $q_1, \ldots, q_{n+1} \in M$ with $p_j \sim q_j$ for $j = 1, \ldots, n+1$ and the construction of the sequence follows by induction.

    Since $\pi$ is normal on $M$, we have that \[\pi\left(\sum_kq_k\right) = \sum_k \pi(q_k).\] 
    Since $\A$ is an AW*-algebra, Lemma \ref{ProjEquivHom} tells us that we must also have \[\pi\left(\sum_k p_k\right) = \sum_k \pi(p_k).\]
    In other words, $\pi$ is CAP on $N$.
    As we know that $\pi$ is faithful, we get that $\pi(N)$ is a von Neumann algebra on $H$.

    Since $N$ was arbitrary, it follows that every MASA of $\pi(\mathcal{A})$ is a von Neumann algebra on $H$.
    By Pedersen's Theorem \ref{PTheorem}, we conclude that $\pi({\mathcal{A}})$ is a von Neumann algebra on $H$.

    It follows that $\mathcal{A} \cong \pi(\mathcal{A})$ is a W*-algebra, as required.
\end{proof} 

\begin{rem}
    Notice that by comparing Theorem \ref{minimalPedersenII1} and Pedersen's Theorem \ref{PTheorem} we see that a dichotomy arises for $\ast$-representations of Type II$_1$ AW*-factors. 
    Namely, for a Type II$_1$ AW*-factor $\mathcal{A}$ and a non-zero $\ast$-representation $(\pi, H)$, either 
    \begin{enumerate}[(i)]
        \item $\pi(M)$ is a concrete von Neumann algebra on $H$ for \emph{all} MASAs $M$ of $\A$, or
        \item $\pi(M)$ is not a concrete von Neumann algebra on $H$ for \emph{any} MASA $M$ of $\A$. 
    \end{enumerate}
\end{rem}

It is not hard to see that we may rephrase our result in the following way, after using Corollary \ref{cor:CAPrep}. 
This rephrasing is perhaps more obviously a strengthened version of Pedersen's Theorem \ref{PTheorem}.

\begin{cor}\label{cor:strongPedII1}
    Let $\A \subseteq \mathcal{B(H)}$ be a Type II$_1$ AW*-factor. Then $\A$ is a von Neumann algebra on $H$ if and only if there exists a MASA of $\A$ which is a von Neumann algebra on $H$.
\end{cor}

\begin{warn}
    Given that the restriction of the quasitrace $\tau$ to a MASA $M$, which we shall denote by $\tau_M$, is a normal, faithful state, it is tempting to try to use a GNS-type argument to construct a representation of $\A$ that is normal on $M$. 
    
    In particular, by Hahn-Banach, we can extend $\tau_M$ to a state $\Phi$ on $\mathcal{A}$. 
    Let $(\pi_\Phi, H)$ be the GNS representation of $\A$ induced by $\Phi$ with the cyclic vector $\xi$, so that $\braket{\pi_\Phi(x)\xi}{\xi}=\Phi(x)$.
    If $x \in M$, then by construction, we have $\Phi(x) = \tau_M(x)$. 
    Since $\tau_M$ is CAP on $M$, it seems that our construction forces $\pi_\Phi$ to be CAP on $M$.
    If this were true then we could argue, as in Theorem \ref{minimalPedersenII1}, that $(\pi_\Phi, H)$ is a representation of $\A$ as a concrete von Neumann algebra on $H$.
    This would seem to show that every Type II$_1$ AW*-factor is a W*-factor (thereby also positively resolving all the equivalent conjectures above). 
    However, this argument is subtly wrong.
    
    While it \emph{is} true that if $(\pi_\tau, K)$ is the GNS representation of $M$ induced by $\tau_M$ then $\pi_\tau$ is CAP on $M$ and $\pi_\tau(M)$ is a von Neumann algebra on $K$, it is possible to extend $\tau_M$ to a state $\Phi$ on $\A$ in a degenerate way, so that the representation $(\pi_\Phi, H)$ is not normal when restricted back to $M$.
    Indeed, even in the case where $\A$ is already known to be a W*-algebra, the degenerate case can arise. 
    Akemann and Sherman showed (see \cite[Lemma 4.3]{AkemannSherman}) that for any Type II$_1$ W*-algebra $\A$ and a MASA M, $\tau_M$ can always be extended to a \emph{singular} state on $\A$. 
    (If $\A$ is a W*-algebra, we can write $\mathcal{A}^* = \mathcal{A}_* \oplus \mathcal{A}_*^{\perp}$. The functionals belonging to $\mathcal{A}_*$ are the normal ones and those belonging to $\mathcal{A}_*^{\perp}$ are the singular ones.)
    Clearly, the GNS representation of any W*-algebra induced by a singular state does not give rise to a normal representation of the W*-algebra as a concrete von Neumann algebra.
    \end{warn}

\section{{Properly Infinite Case}}\label{ProperlyInfiniteSec}
Finally, we shall prove some partial results pointing towards the W*-Pedersen characterization in the case of properly infinite AW*-algebras. 
We require a Mackey-Gleason type result for certain unital C*-algebras in order to proceed. 

\begin{defn}
       Let $\A$ be a unital C*-algebra. A \emph{quasilinear functional} on $\A$ is a function $\rho: \mathcal{A} \to \mathbb{C}$ such that 
    \begin{enumerate}[(i)]
        \item $\rho$ is bounded on the unit ball of $\mathcal{A}$,
        \item $\rho$ is a linear functional on every commutative C*-subalgebra of $\mathcal{A}$,
        \item $\rho(a+bi) = \rho(a) + i\rho(b)$ for any $a, b \in \mathcal{A}_{sa}$.
    \end{enumerate}
    We say that a quasilinear functional $\rho$ is 
    \begin{enumerate}[(i)]
        \item \emph{continuous} if $\rho(x_i) \to \rho(x)$ whenever $(x_i) \to x$ in norm, 
        \item \emph{positive} if the restriction of $\rho$ to each commutative C*-subalgebra is a positive linear functional, 
        \item \emph{completely additive on projections} (CAP) if it is both positive and preserves suprema of arbitrary families of mutually orthogonal projections.
    \end{enumerate}
    Finally, we call $\rho$ a \emph{quasistate} if $\rho(1) = 1$. 
\end{defn}

\begin{prop}\cite[Theorem 4.1]{Kaplan} \label{Kaplan} 
Let $\mathcal{A}$ be a unital C*-algebra and let $\rho$ be a positive quasistate on $\A$ that is linear when restricted to each C*-subalgebra generated by two projections. Suppose that there exists a sequence of projections $(p_n)_{n=1}^{\infty} \subseteq \Proj(\mathcal{A})$ such that
\begin{enumerate}
    \item $\rho(1-p_n) \to 0$, and 
    \item for each $n$ there are partial isometries $u_n, v_n,$ and $w_n$ in $\A$ such that $u_n^*u_n = v_n^*v_n = w^*_nw_n = p_n$ and $p_n, u_nu_n^*, v_nv_n^*, w_nw_n^*$ are mutually orthogonal.
\end{enumerate}
    Then there is a state on $\A$ that coincides with $\rho$ on $\Proj (\mathcal{A})$. Moreover, if $\A$ contains a dense set of elements with finite spectrum and $\rho$ is continuous, then $\rho$ is linear.
\end{prop}

We can now state and prove a Mackey-Gleason result for properly infinite AW*-algebras. 

\begin{thm}\label{thm:properlyInfiniteMackeyGleason}
     Let $\mathcal{A}$ be a properly infinite AW*-algebra. If $\A$ admits a CAP quasistate then $\A$ admits a CAP state. Moreover, any continuous CAP quasistate on $\A$ is linear.
\end{thm}

\begin{proof}
    We proceed by checking that the conditions of Proposition \ref{Kaplan} are met.
    
    If $\mathcal{A}$ is properly infinite, then by \cite[Lemma 4.4]{Kap51} there exists an orthogonal sequence of projections $(e_n)_{n=1}^\infty$ such that $\sup \{e_n\} = 1$ and $e_n \sim 1$ for all $n \geq 1$. 
    Define $p_n := e_1 + \cdots + e_n$ for each $n \geq 1$.
    Then $p_n \sim 1$ for each $n \geq 1$ and if we define $q_n = e_{n+1}$, $s_n := e_{n+2}$, and $t_n := e_{n+3}$, we can find partial isometries in $\A$ as in the statement of Proposition \ref{Kaplan} (2) from equivalence of the projections and the fact that $p_n, q_n, s_n$ and $t_n$ are mutually orthogonal.

    Let $\rho: \mathcal{A} \to \mathbb{C}$ be a CAP quasistate on $\mathcal{A}$. 
    Define $r_1 := p_1$ and $r_k := p_k - p_{k-1} = p_k(1-p_{k-1})$. 
    By construction $p_{k-1} \leq p_k$, so each $r_k$ is a projection, the $r_k$'s are mutually orthogonal, with \[p_n = \sum_{k=1}^n r_k \qquad \text{ and } \qquad \sum_{k \geq 1}r_k = \bigvee_{n \geq 1}p_n = 1.\]
    So, for every $n$, we have $\displaystyle 1-p_n = \sum\limits_{k > n} r_k$.
    Since $\rho$ is CAP, 
    \begin{align*}
        \rho(1) = \rho\left(\sum_{k=1}^{\infty} r_k\right) = &\sum_{k=1}^{\infty} \rho(r_k) \\
        \implies &\sum\limits_{k > n} \rho(r_k) \to 0.
    \end{align*}
    In particular, $\rho(1-p_n) \to 0$ as $n \to \infty$ and our sequence also satisfies condition $(1)$ of Proposition \ref{Kaplan}.
    
    Hamhalter showed \cite[Corollary 3.9]{Hamhalter} that in any AW*-algebra $\B$ without Type I$_2$ direct summand, if $\rho$ is any quasilinear functional on $\B$ and $p, q \in \Proj(\mathcal{B})$, then $\rho$ is linear on the AW*-algebra generated by $p, q$. 
    That is, $\rho$ is linear on $AW^*(p,q) \subseteq \mathcal{B}$. 
    Since our AW*-algebra $\A$ is properly infinite, it has no Type I$_2$ direct summand and so $\rho$ is linear on $AW^*(p,q)$ for any $p, q \in \Proj(\mathcal{A})$.
    Finally, $C^*(p,q) \subseteq AW^*(p,q)$, so $\rho$ must be linear on any C*-algebra generated by two projections in $\mathcal{A}$.
    By Proposition \ref{Kaplan}, there exists a state $\Phi: \mathcal{A} \to \mathbb{C}$ such that $\Phi |_{\Proj(\mathcal{A})} = \rho$. In other words, $\Phi$ is a CAP state.  
    
    If $\rho$ is continuous, then $\Phi = \rho$ follows from the ``moreover'' part of the statement of Proposition \ref{Kaplan}. 
    As the details in the reference are somewhat sparse, we complete the argument below. 

    If $a \in \mathcal{A}_{sa}$ has finite spectrum, then taking its spectral decomposition in $C^{*}(1, a)$, we get $a = \sum_{k=1}^n \lambda_k p_k$ for some $p_k \in \Proj(\mathcal{A})$ and $\lambda_k \in \mathbb{R}$.
    Since $\rho$ is linear on $C^*(1,a)$, we get \[\rho(a) = \sum_{k=1}^n \lambda_k \rho(p_k) = \sum_{k=1}^n \lambda_k \Phi(p_k) = \Phi(a).\]

    Since AW*-algebras have RR0, any self-adjoint element of $\A$ can be approximated in norm by elements of finite spectrum. 
    By continuity of $\rho$, and the agreement of $\rho$ and $\Phi$ on self-adjoint elements of finite spectrum, we deduce that $\rho(a) = \Phi(a)$ for any $a \in \mathcal{A}_{sa}$.
    Writing an arbitrary $x \in \mathcal{A}$ as $x = a+ib$ where $a, b \in \mathcal{A}_{sa}$, we see that \[\rho(x) = \rho(a+ib) = \rho(a) +i\rho(b) = \Phi(a)+i\Phi(b) = \Phi(x).\]
    Therefore, $\rho = \Phi$ as desired.
\end{proof}

\begin{rem}
    During the preparation of this manuscript, the author came across a closely-related result in a Russian-language paper by Matvejchuk \cite{Mat86}. 
    According to the MathSciNet review, which is written in English, the main result of that paper is that any \emph{finitely additive measure} on a purely infinite AW*-algebra can be extended to a linear functional on the algebra.
    There, a finitely additive measure on an AW*-algebra $\A$ is a map $\mu: \Proj(\A) \to \mathbb{R}_+$ satisfying  \[ \mu\left(\sum_{i} p_i\right) = \sum_{i} \mu(p_i),\] where $(p_i)$ is a finite family of mutually orthogonal projections. 
    It is possible that Theorem \ref{thm:properlyInfiniteMackeyGleason} follows in a straightforward way from that result, depending on what is meant by ``purely infinite'' in \cite{Mat86} (cf. Remark \ref{rem:infiniteTerminology}). 
    
    Because the text is only available in Russian, we were unable to verify the exact hypotheses imposed on the AW*-algebra or check the details of the proof.
    Unfortunately, this leaves us unable to say more about the exact relationship between the preceding result and the one given in \cite{Mat86}.
\end{rem}

\begin{cor}
    Let $\A$ be a properly infinite AW*-algebra. If $\A$ admits a separating family of CAP quasistates, then $\A$ is a W*-algebra.
\end{cor}

\begin{proof}
    The given family of CAP quasistates $(\rho_i)_{i \in I}$ gives rise to a family of states $(\Phi_i)_{i \in I}$ agreeing with $\rho_i$ on $\Proj(\A)$ by Theorem \ref{thm:properlyInfiniteMackeyGleason}.
    It remains to show that separation is preserved on arbitrary positive elements.

    For any $0 \neq a \in \A_+$, $\varepsilon \in (0, \norm{a})$ we can find $0 \neq p \in \Proj(\A)$ with $a \geq \varepsilon p,$ by Corollary \ref{cor:dominateProjection}.
    Because the quasistates separate $\A_+,$ there is some $\rho_i$ such that $\rho_i(p) > 0.$
    Therefore 
    \[ \Phi_i(a) \geq \varepsilon \Phi_i(p) = \varepsilon \rho_i(p) > 0.\]
    
    So the family of CAP states $\Phi_i$ is separating and the fact that $\A$ is a W*-algebra follows by Corollary \ref{CAPcor}.
\end{proof}

By the above results, we would like to try and show something like ``if $\A$ is properly infinite and every MASA is a W*-algebra, then $\A$ admits a separating family of CAP quasistates'', from which we could conclude that $\A$ is a W*-algebra. 

Fixing a MASA $M_i$ of such an $\mathcal{A}$, we can certainly take a separating family of normal (on $M_i$) states \[(\varphi^{(i)}_\lambda)_{\lambda \in \Lambda} \subseteq (M_i)_*\] and, by Hahn-Banach, extend them to a family of states \[(\Phi^{(i)}_\lambda)_{\lambda \in \Lambda} \subseteq \mathcal{A}^*.\] 
Repeating this process for all MASAs, we get a separating family of states \[(\Phi^{(i)}_\lambda)_{i, \lambda} \subseteq \A^*\] such that for any $\lambda \in \Lambda$, each $i$-indexed state $\Phi^{(i)}_\lambda$ was obtained as Hahn-Banach extension of a normal state on $M_i$.
However, there is no \emph{a priori} reason why every state $\Phi^{(i)}_\lambda$ will be normal on every MASA, so we have not exhibited a separating family of CAP states. 
While the restriction of $\Phi_\lambda^{(i)}$ back to $M_i$ will be normal, there is no \emph{a priori} reason why $\Phi_\lambda^{(i)}$ should be CAP on all other MASAs. 
It is possible that (a subfamily of) this family of states can be chosen -- or perhaps perturbed -- to find a separating family of CAP quasistates, so that the above results can be exploited, but it is not clear to the author how to proceed along these lines.

Nevertheless, if we restrict to the case of countably decomposable properly infinite AW*-factors, we can at least obtain a strengthening of Pedersen's Theorem akin to Theorem \ref{minimalPedersenII1}, by taking a different approach.
To proceed, we will need the following lemma.

\begin{lem}\label{equivCountableProj}
    Let $\A$ be a countably decomposable, properly infinite AW*-algebra and let $p, q \in \Proj(\mathcal{A})$ be properly infinite projections. If $z(p) \leq z(q)$ then $p \lesssim q$. In particular, in a countably decomposable, properly infinite AW*-factor, all properly infinite projections are equivalent. 
\end{lem}

\noindent A standard proof of this result in the von Neumann algebra case relies only on properties of projections shared by both (properly infinite, countably decomposable) AW*-algebras and W*-algebras.
In particular, the only properties used are the following:
\begin{enumerate}[(i)]
    \item \cite[Proposition 12.1]{Berberian} the generalized Cantor-Schr{\"o}der-Bernstein theorem holds for projections;
    \item \cite[Corollary 14.1]{Berberian} any AW*-algebra $\A$ has \emph{generalized comparison}  meaning that for any two projections $e, f$ in $\A$ there exists a central projection $h$ such that $he \lesssim hf$ and $(1-h)f \lesssim (1-h)e$; and 
    \item \cite[Theorem 17.1]{Berberian} in a properly infinite AW*-algebra, there exists some non-zero projection $g$ such that $g \sim 1-g \sim 1.$ 
\end{enumerate}
Therefore, the proof follows verbatim as in \cite[Proposition 3.2.9]{Peterson}, so we omit the details here.

\subsection{Simple Factors} In the case of countably decomposable, simple C*-algebras which are AW*-factors, we get a particularly nice Pedersen-type characterization of the W*-factors, see Theorem \ref{mainTheorem2}. The non-simple case will be treated separately.

\smallskip
We first consider the countably decomposable Type III AW*-factors; these are particularly interesting because the first known examples of  \emph{wild} AW*-factors (i.e., AW*-factors that are not W*-factors) were of this kind \cite{SaitoFactors}.
Moreover, Hamana has shown that there exist many such examples. In fact, he showed that there are $2^{2^{\aleph_0}}$ non-isomorphic countably decomposable wild Type III AW*-factors \cite{Hamana}.

\begin{prop}\label{TypeIIIFactor}
    Let $\A$ be a countably decomposable Type III AW*-factor. Then $\A$ is a W*-factor if and only if 
    \begin{enumerate}[(i)]
        \item every MASA of $\A$ is a W*-algebra, and 
        \item there exists a $\ast$-representation $(\pi, H)$ of $\A$ and a MASA $M$ of $\mathcal{A}$ such that $\pi|_M$ is non-zero and CAP.
    \end{enumerate}
\end{prop}

\begin{proof}
    The forward direction is trivial; we show only the converse.
    
    First, since $\A$ is a simple C*-algebra, we must have that either $\pi$ is the zero $\ast$-homomorphism or it is faithful.
    By our assumption that $\pi|_M$ is non-zero, we must have that $\pi$ is a faithful $\ast$-homomorphism of $\mathcal{A}$ into $\mathcal{B}(H)$.

    Now, suppose that $\A$ is a countably decomposable Type III AW*-factor, that every MASA of $\A$ is a W*-algebra, and that $(\pi, H)$ is a faithful $\ast$-representation of $\A$ such that $\pi|_M$ is non-zero and CAP.
    Then $\pi(M) \cong M$ is a von Neumann algebra on $H$.
    
    Let $N$ be an arbitrary MASA of $\mathcal{A}$. 
    Let $(q_k)$ be a family of orthogonal projections in $N$ with supremum $q$.
    By countable decomposability, we may assume that the family is countable, hence indexed by $k \in \mathbb{N}$.
    Let $S := \{ k \in \mathbb{N}: q_k \neq 0\}.$
    
    For each $q_k$ with $q_k = 0$ (i.e., with $k \in \mathbb{N} \setminus S$), set $p_k := 0 \in M$.
    Since $q_k = p_k = 0$ for these indices, we have $q_k \sim p_k$ for all $k \in \mathbb{N} \setminus S.$
    If $S = \emptyset$, i.e., if $(q_k)_{k \in \mathbb{N}}$ is a sequence of all $0$'s, we are done.
    If not, we still need to construct mutually orthogonal projections $p_k$ with $q_k \sim p_k$ for $k \in S$. 

    By \cite[Theorem 6.8]{HR1}, for any set $S$ with cardinality $1 \leq |S| \leq \aleph_0$ and a fixed MASA $M$, we can find an orthogonal family of projections $(p_k)_{k \in S} \subseteq \Proj(M)$  such that $$\sup_{\Proj(\A)} \{p_k\} = 1$$ and $p_k \sim 1$ for all $k$.
    
    Recall that in any countably decomposable Type III AW*-factor, every non-zero projection is equivalent (see Lemma \ref{equivCountableProj}).
    This means that the family of orthogonal projections $(p_k)_{k \in S}$ in $M$ satisfies $q_k \sim p_k$ for all $k \in S$.
    By the preceding construction, this means that we have $q_k \sim p_k$ for all $k \in \mathbb{N}.$

    By assumption $(ii)$, \[\pi\left(\sup \{p_k\}\right) = \sum_k \pi(p_k),\] where the right-hand sum refers to the SOT-sum in $\B(H)$.  
    By Lemma \ref{ProjEquivHom}, we also have \[\pi(q) = \pi\left(\sup\{q_k\} \right) = \sum_k \pi(q_k).\] 
    So $\pi |_N$ is CAP, hence $\pi(N)$ is a von Neumann algebra on $H$.
    By Pedersen's Theorem \ref{PTheorem}, $\pi(A)$ must be a von Neumann algebra on $H$. 
    Since $\pi$ is faithful, it follows that $\mathcal{A} \cong \pi(\mathcal{A}) \subseteq \mathcal{B}(H)$, where the latter is a concrete von Neumann algebra on $H$, so we are done. 
\end{proof}

\begin{rem}\label{rem:masaCondition}
    Condition $(i)$ in Proposition \ref{TypeIIIFactor} is superfluous. 
    As long as we have a MASA $M$ and a faithful, CAP representation of $M$ on $H$, then $M$ is a von Neumann algebra on $H$ (by Corollary \ref{cor:CAPrep}, since $M$ is necessarily an AW*-algebra under our hypotheses).
    We have left condition $(i)$ as written above, to make the parallel to the concrete Pedersen Theorem \ref{PTheorem} more obvious.
\end{rem}

As in the Type II$_1$ AW*-factor case, we can rephrase our result as follows. 

\begin{cor}\label{cor:strongPedIII}
    Let $\A \subseteq B(H)$ be a countably decomposable Type III AW*-factor. Then $\A$ is a von Neumann algebra on $H$ if and only if there exists a MASA of $\A$ which is a von Neumann algebra on $H$.
\end{cor}

We may now state and prove our second main theorem.

\begin{thm}\label{mainTheorem2}
     Let $\A$ be a simple unital C*-algebra in which every family of mutually orthogonal projections is at most countably infinite.  Then $\A$ is a W*-factor if and only if \begin{enumerate}[(i)]
        \item every MASA of $\A$ is an AW*-algebra, and 
        \item there exists a $\ast$-representation $(\pi, H)$ of $\A$ and a MASA $M$ of $\mathcal{A}$ such that $\pi|_M$ is non-zero and CAP.
    \end{enumerate}
\end{thm}

\begin{proof}
     The forward direction is clear, since every simple W*-factor is either Type $I_n$ for some $n \in \mathbb{N}$, Type II$_1$, or countably decomposable of Type III. 
     
     For the converse, notice that condition $(i)$ says that $\A$ is an AW*-algebra and simplicity says that $\A$ has trivial centre. 
     Recalling that the only countably decomposable, simple AW*-factors are those of Type I$_{n}$ for some $n \in \mathbb{N}$, Type II$_1$, or countably decomposable of Type III, the result follows by collecting our previous results.

      In particular, we know that \emph{every} Type I factor is a W*-algebra by \cite[Theorem 2]{Kap52} (cf. Proposition \ref{TypeI}), the Type II$_1$ case is handled by Theorem \ref{minimalPedersenII1}, and the countably decomposable Type III case was shown in Proposition \ref{TypeIIIFactor}.
\end{proof}

\begin{cor}
    Let $\A \subseteq B(H)$ be a countably decomposable, simple AW*-factor. Then $\A$ is a von Neumann algebra on $H$ if and only if there exists a MASA of $A$ which is a von Neumann algebra on $H$.
\end{cor}

\subsection{Non-simple Factors} The case of Type II$_\infty$ AW*-factors is a bit more delicate, but we can still obtain a partial strengthening of Pedersen's Theorem. In this case, we are able to remove the assumption of countable decomposability.
\smallskip

For this case, we will need to use the fact that any AW*-factor $\A$ is \emph{normal} \cite[Corollary 4.7]{SWnormal}. 
Essentially, this means that the projection lattice of $\A$ sits ``nicely'' inside $\Asa$.

\begin{defn}
    Let $\A$ be an AW*-algebra. We say that $\A$ is \emph{normal} if for any directed set of projections $(p_i)_{i \in I}$ in $\A$ with supremum $p \in \Proj(\mathcal{A}),$ the supremum of $(p_i)$ in $\Asa$ exists and is equal to $p$. That is, if $a \in \mathcal{A}_{sa}$ and $a \geq p_i$ for all $i \in I$, then $a \geq p$.
\end{defn}

 In addition to AW*-factors, it is known that every monotone complete C*-algebra is a normal AW*-algebra. However, the question of whether arbitrary AW*-algebras are normal is open.
\smallskip 

Finally, we will need the following result, which holds in arbitrary AW*-algebras.

\begin{lem}\cite[Proposition 2.1.10]{SWmonograph}\label{lem:conjLim}
    Let $\A$ be an AW*-algebra. If $(a_i)_{i \in I}$ is a norm-bounded, monotone increasing net of elements in $\Asa$ with supremum $a \in \Asa$ and $x \in \A$, then $(x^*a_ix)_{i \in I}$ has supremum $x^*ax \in \Asa.$
\end{lem}

Note that the statement assumes that we are given a norm-bounded monotone increasing net of self-adjoint elements with supremum, \emph{not} that the supremum of any norm-bounded monotone increasing net of self-adjoint elements exists. 

\smallskip
We are now able to show our strengthened Pedersen Theorem for Type II$_\infty$ AW*-factors.
In the proof, we use our usual $\sup_{\Asa}\{ \cdot \}$ notation for suprema of norm-bounded, increasing families of self-adjoint elements and reserve the $\vee (\cdot)$ notation for suprema of projections, to avoid confusion.
\smallskip

\begin{prop}\label{TypeIIinfinityfactor}
    Let $\A$ be a Type II$_\infty$ AW*-factor. Then $\A$ is a W*-factor if and only if \begin{enumerate}[(i)]
        \item every MASA of $\A$ is a W*-algebra, and 
        \item there exists a MASA $M$ which contains a non-zero finite projection $e$ and a $\ast$-representation $(\pi, H)$ of $\A$ such that $\pi|_M$ is non-zero, CAP, and $\pi(e) \neq 0$.
    \end{enumerate}
\end{prop}

\begin{proof}
    The forward direction is immediate.

    Let $\A$ be a Type II$_\infty$ AW*-factor and let $(\pi, H)$, $M \subseteq A$, and $e$ be as in the assumption.
     Since $e$ is finite, $e\A e$ is a Type II$_1$ AW*-factor and $eMe$ is a MASA of $e \A e$. 
    
     Let $\pi_e: e\A e \to \mathcal{B}(\pi(e)H)$ be the $\ast$-representation given by \[ \pi_e(a) = \pi(a)|_{\pi(e)H}.\] 
    Since $\pi_e(e) = I|_{\pi(e)H}$ and $\pi(e) \neq 0$, this representation is non-zero.  
    Moreover, the restriction of the representation $\pi_e|_{eMe}$ is CAP, because $\pi|_M$ is CAP. 
    
    In particular, $e \A e$ is a Type II$_1$ AW*-factor that admits a faithful $\ast$-representation $(\pi_e, \pi(e)H)$ with a MASA $eMe$ which is a von Neumann algebra on $\pi(e)H$.
    Therefore, Theorem \ref{minimalPedersenII1} tells us that $e \A e$ is a W*-factor.
    In particular, $e \A e$ admits a normal faithful tracial state $\varphi: e\A e \to \mathbb{C}.$
       
    For each $x \in \mathcal{A}$ this allows us to define $\Phi_x: \A \to \mathbb{C}$ by \[\Phi_x(a) = \varphi\left(ex^*axe\right).\]
    It is clear that each $\Phi_x$ is a bounded positive linear functional with \[\norm{\Phi_x} = \Phi_x(1) = \varphi(ex^*xe).\]
    Moreover, by faithfulness of $\varphi$ we see that $\Phi_x$ is not identically the zero functional, unless $xe = 0.$
    We claim that $(\Phi_x)_{x \in \A}$ is a separating family of CAP linear functionals on $\A$.
    
    First, we show the separation property.
    Towards this goal, let $0 \neq a \in \A_+$ and pick $\varepsilon \in (0, \norm{a})$.
    By Corollary \ref{cor:dominateProjection}, there is some $0\neq p \in \Proj(\A)$ such that $a \geq \varepsilon p$ and by generalized comparability in a factor, we have either \[ e \lesssim p \quad \text{ or } \quad p \lesssim e.\] 
    In either case, there exist non-zero projections $q, f \in \Proj(\A)$ such that $q \leq p$ and $ f \leq e, $ and there is a partial isometry $v \in \mathcal{A}$ such that $q = vv^* \sim v^*v = f.$
    This means that $pv = v$ and $ve = v,$ so 
    \begin{align*}
        \Phi_v(a) = \varphi(ev^*ave) \geq \varepsilon\varphi(ev^*pve) 
        = \varepsilon \varphi(v^*pv) 
        = \varepsilon \varphi(v^*v) 
        = \varepsilon \varphi(f) > 0, 
    \end{align*}
    where the final inequality comes from the fact that $\varphi$ is faithful and $f \neq 0$.
    So $(\Phi_x)_{x \in \A}$ is a separating family of bounded linear functionals.
    
    It remains to show that each functional is CAP.
    For this, let us consider $(p_i)_{i \in I}$, a family of mutually orthogonal projections with
   \[p = \bigvee_{i \in I} p_i = \sup_{F \in \mathcal{F}}\sum_{i \in F} p_i,\]
   where $\mathcal{F}$ denotes the directed family of finite subsets of $I$. 
   For each $F \in \mathcal{F}$, write \[s_F := \sum_{i \in F} p_i,\] and notice that $(s_F)_{F \in \mathcal{F}}$ is a norm-bounded, monotone increasing net in $\Asa.$
   Since $\A$ is a normal AW*-algebra, it follows that 
   \[ p = \sup_{\Asa} \{s_F: F \in \mathcal{F}\}.\]
   Moreover, for any fixed $x \in \A$, Lemma \ref{lem:conjLim} gives \[\sup_{\Asa} \{(xe)^*s_F(xe): F \in \mathcal{F}\} = (xe)^*p(xe) = ex^*pxe.\] 
   Notice that this is also the supremum of $\left(ex^*s_Fxe\right)_{F \in \mathcal{F}}$ as a net in $(e \A e)_{sa}$.
    
   Putting everything together, we have for any fixed $x \in \A$,
    \begin{align*}
        \Phi_x(p) = \varphi(ex^*pxe) = \varphi\left(\sup_{(e \A e)_{sa} }\left\{(xe)^*s_F(xe): F \in \mathcal{F}\right\}\right) &= \sup_F \varphi(ex^*s_Fxe) \\
        &= \sup_F \varphi\left(ex^*\left(\sum_{i \in F} p_i \right)xe\right) \\
        &= \sup_{F} \sum_{i \in F} \varphi(ex^*p_ixe) \\
        &= \sup_{F} \sum_{i \in F} \Phi_x(p_i).
    \end{align*}
    This means that \[ \Phi_x \left(\bigvee_{i \in I}p_i\right) = \sum_{i \in I} \Phi_x(p_i),\]
    hence $\Phi_x$ is CAP.
    
    After normalizing the non-zero functionals, i.e., the functionals $(\Phi_x)_{x \in \A}$ where $xe \neq 0$, we are left with a separating family of CAP states on $\A$, so $\A$ is a W*-algebra by Corollary \ref{CAPcor}.
\end{proof}

\begin{rems} Two remarks are in order. 
    \begin{enumerate}[(i)]
        \item Proposition \ref{TypeIIinfinityfactor} does not obviously prove that the representation $(\pi, H)$ of $\A$ given in condition $(ii)$ is a concrete faithful representation of $\A$ as a von Neumann algebra on $H$. We only establish the abstract fact that $\A$ must be a W*-algebra. As such, it is not clear that this result can be rephrased in a more concrete way, as we did for countably decomposable, simple AW*-factors. 
        \item Again, condition $(i)$ in Proposition \ref{TypeIIinfinityfactor} is not necessary. We essentially only needed condition $(ii)$ to show that $e \A e$ is a W*-algebra, then we proved that this implies that $\A$ is a W*-algebra. We left the first condition in the statement above, in order to draw parallels with the concrete Pedersen Theorem \ref{PTheorem}; we remove it below for full generality.
    \end{enumerate}
\end{rems}

For our third and final main result, we return to the assumption that $\A$ is countably decomposable. 
We do this to exclude the case of arbitrary Type III AW*-factors, which we have not studied in this paper.

\begin{thm}\label{mainThm3}
    Let $\A$ be a unital, non-simple C*-algebra with trivial centre and in which every family of mutually orthogonal projections is at most countably infinite. Then $\A$ is a W*-factor if and only if
    \begin{enumerate}[(i)]
        \item every MASA of $\A$ is an AW*-algebra, and 
        \item there exists a MASA $M$ which contains a non-zero finite projection $e$ and a $\ast$-representation $(\pi, H)$ of $\A$ such that $\pi|_M$ is non-zero, CAP, and $\pi(e) \neq 0$.
    \end{enumerate}
\end{thm}

\begin{proof}
    As noted in Section \ref{subsection:AW}, the only non-simple, countably decomposable AW*-factors are those of Type I$_\infty$ and Type II$_\infty$.
    The converse is then immediate from the fact that all Type I AW*-factors are W*-algebras (\cite[Theorem 2]{Kap52}) and Proposition \ref{TypeIIinfinityfactor}.

    The forward direction follows from the fact that any non-simple W*-factor is of Type I$_\infty$ or Type II$_\infty$ and therefore contains a finite projection. 
\end{proof}

Before we conclude, let us return to the discussion on C*-quotients of W*-algebras that was initiated in Section \ref{bidualChar}.

\section{Monotone Completeness}\label{MCQuestion}

This section was inspired by a MathOverflow discussion between Hannes Thiel and David Gao \cite{MO}. 
Therein, the former gave our previously established Example \ref{bidualRetractEx} and a sketch of a proof for Proposition \ref{bidualRetractProp}; the latter pointed out a relevant commutative example in the literature, which we include here as Example \ref{exmp:injectiveAbelian}.

\smallskip
Recall that in Section \ref{bidualChar}, we claimed that there are AW*-algebras which can be realized as the C*-quotient of a W*-algebra but are not themselves W*-algebras. We consider this example now.

\begin{exmp}\label{exmp:injectiveAbelian}
     Let $\mathcal{C}(X)$ be any commutative AW*-algebra with $X$ Stonean but \emph{not} hyperstonean, so that $\mathcal{C}(X)$ is a non-W*, AW*-algebra (as in Example \ref{exmp:nonWseparable}). 
    Then $\mathcal{C}(X)$ is injective in the category of unital, commutative C*-algebras with $\ast$-homomorphisms \cite[Theorem 2.4]{HadwinPaulsen}.
    In this case, viewing $\mathcal{C}(X)$ as its canonical copy inside $\mathcal{C}(X)^{**}$, injectivity means that the identity map $\iota:\mathcal{C}(X) \to \mathcal{C}(X)$ can be extended to a $\ast$-homomorphism $\hat{\iota}: \mathcal{C}(X)^{**} \twoheadrightarrow \mathcal{C}(X)$.
    Then $\mathcal{C}(X) \cong \mathcal{C}(X)^{**}/{\ker\hat{\iota}}$, but we know that $\mathcal{C}(X)$ is not a W*-algebra. 
\end{exmp}
   
The key idea here is that $\mathcal{C}(X)$ sits nicely inside its bidual $\mathcal{C}(X)^{**}$. 
In particular, $\mathcal{C}(X)$ is a $\ast$-homomorphic retract of its bidual; this example motivated Definition \ref{defn:retract} introduced above.
 \smallskip

As Example \ref{exmp:injectiveAbelian} highlights, a commutative unital C*-algebra being a retract of its bidual does not always imply that it is a W*-algebra (cf. the Type II$_1$ AW*-factor case in Corollary \ref{cor:bidualII_1}).
Recalling that every commutative AW*-algebra is monotone complete, we come to a result that genuinely extends the commutative example.

\begin{prop}\label{bidualRetractProp}
    If a unital C*-algebra is a retract of its bidual, then it is monotone complete.
\end{prop}

\begin{proof}
    Suppose that $\A$ is a unital C*-algebra for which there exists a multiplicative conditional expectation $E: \mathcal{A}^{**} \twoheadrightarrow \mathcal{A}$. 
    Let $(x_i)_{i \in I}$ be a norm-bounded, monotone increasing net in $\mathcal{A}_{sa}$.
    Since $\mathcal{A}^{**}$ is monotone complete, there is some $y \in (\mathcal{A}^{**})_{sa}$ such that \[y = \sup_{(\A^{**})_{sa}}\{x_i: i \in I\}.\]
    Since $E$ preserves order and fixes $\mathcal{A}$, we have $x_i = E(x_i) \leq E(y)$, so $E(y)$ is an upper bound for $(x_i)$.
    We claim that \[E(y) = \sup_{\A_{sa}} \{x_i: i \in I\}.\] 
    
    Let $x \in \mathcal{A}_{sa}$ be any upper bound for $(x_i)$ in $\mathcal{A}_{sa}$.
    Then since $y$ is the LUB in $(\A^{**})_{sa}$, we have $y \leq x$.
    By order preservation, we have $E(y) \leq E(x) = x$. 
    So $(x_i)$ has a LUB, namely $E(y)$ in $\A_{sa}$.
\end{proof}

By Example \ref{exmp:injectiveAbelian}, we know that the converse of Proposition \ref{bidualRetractProp} holds in the commutative case. This raises the following question.

\begin{quest}\label{q:Bidual}
    Is every monotone complete C*-algebra a retract of its bidual?
\end{quest}

A positive answer would mean that every monotone complete Type II$_1$ AW*-factor is a W*-factor (see the proof of Corollary \ref{cor:bidualII_1}).
In this case, the resolution of Kaplansky's Conjecture, the 2Q Conjecture, and the finite W*-Pedersen Conjecture would be reduced to the following.

\begin{quest}\label{q:MC}
    Is every Type II$_1$ AW*-factor monotone complete?
\end{quest}

Unfortunately, the question of whether arbitrary AW*-algebras are monotone complete is the other longstanding open problem in the theory of AW*-algebras (see e.g., the discussion in \cite[Chapter 8]{SWmonograph}). 
However, recalling the Sait{\^o}-Wright Characterization \ref{AWcharacterization} that a unital C*-algebra is an AW*-algebra if and only if each of its MASAs is monotone complete, the question of monotone completeness of AW*-algebras can be rephrased in a very familiar way (cf. Question \ref{q:commW}).

\begin{quest}
    Is the ``extra'' structure that turns an AW*-algebra into a monotone complete C*-algebra determined commutatively?
\end{quest}

\end{document}